\documentstyle[12pt,amssymb,leqno]{article}
\def\fn#1{\mathop{{\rm #1}\vphantom{\dim}}}

\makeatletter
\renewcommand{\section}{\@startsection{section}{1}{0mm}{8mm}{3mm}{\bf \raggedright}}
\def\@citex[#1]#2{\if@filesw\immediate\write\@auxout{\string\citation{#2}}\fi
  \def\@citea{}\@cite{\@for\@citeb:=#2\do
    {\@citea\def\@citea{\@citesep}\@ifundefined
       {b@\@citeb}{{\bf ?}\@warning
       {Citation `\@citeb' on page \thepage \space undefined}}%
{\csname b@\@citeb\endcsname}}}{#1}}
\def\@citesep{; }
\makeatother

\newtheorem{theorem}{\indent Theorem}[section]

\newtheorem{defn}[theorem]{\indent Definition}
\newtheorem{prop}[theorem]{\indent Proposition}
\newtheorem{lemma}[theorem]{\indent Lemma}
\newtheorem{coro}[theorem]{\indent Corollary}
\newtheorem{remark}{\indent Remark.}[section]

\newtheorem{case}{\indent Case}[section]

\begin{document}

\centerline {\bf RETRACT RATIONALITY AND NOETHER'S PROBLEM}

\

\centerline {Ming-chang Kang}

\centerline {Department of Mathematics}

\centerline {National Taiwan University}

\centerline {Taipei, Taiwan, Rep. of China}

\centerline {E-mail: kang@math.ntu.edu.tw}

\

\vskip 10mm \baselineskip 0.7cm {\bf Abstract.} Let $K$ be any
field and $G$ be a finite group.  Let $G$ act on the rational
function fields $K(x_g:g \in G)$ by $K$-automorphisms defined by
$g \cdot x_h=x_{gh}$ for any $g, \, h \in G$. Denote by $K(G)$ the
fixed field $K(x_g: g \in G)^G$.  Noether's problem asks whether
$K(G)$ is rational (=purely transcendental) over $K$.  We will
prove that, if $K$ is any field, $p$ an odd prime number, and $G$
is a non-abelian group of exponent $p$ with $|G|=p^3$ or $p^4$
satisfying $[K(\zeta_p):K] \leq 2$, then $K(G)$ is rational over
$K$. A notion of retract rationality is introduced by Saltman in
case $K(G)$ is not rational. We will also show that $K(G)$ is
retract rational if $G$ belongs to a much larger class of
$p$-groups. In particular, generic $G$-polynomials of $G$-Galois
extensions exist for these groups.

\vskip 2cm

\noindent $\underline {\mskip 200mu}$



\newpage
\section{Introduction}

Let $K$ be any field and $G$ be a finite group.
Let $G$ act on the rational function field $K(x_g:g\in G)$ by $K$-automorphisms defined by $g\cdot x_h=x_{gh}$ for any $g,h\in G$.
Denote $K(G):=K(x_g:g\in G)^G$ the fixed subfield under the action of $G$.
Noether's problem asks whether $K(G)$ is rational (= purely transcendental) over $K$.

Noether's problem for abelian groups was studied by Fischer,
Furtw\"{a}ngler, K.\ Masuda, Swan,\ Voskresenskii, S.\ Endo and
T.\ Miyata, Lenstra, etc. It was known that $K({\Bbb Z}_p)$ is
rational if ${\Bbb Z}_p$ is a cyclic group of order $p$ with
$p=3,5,7$ or 11. The first counter-example was found by Swan:
${\Bbb Q}({\Bbb Z}_{47})$ is not rational over ${\Bbb Q}$
\cite{Sw1}. However Saltman showed that ${\Bbb Q}({\Bbb Z}_p)$ is
retract rational over $K$ for any prime number $p$, which is
enough to ensure the existence of a generic Galois $G$-extension
and will fulfill the original purpose of Emmy Noether \cite{Sa1}.
For the convenience of the reader, we recall the definition of
retract rationality.

\begin{defn}[\cite{Sa3}] \rm
Let $K\subset L$ be a field extension. We say that $L$ is retract
rational over $K$, if there is a $K$-algebra $R$ contained in $L$
such that (i) $L$ is the quotient field of $R$, and (ii) the
identity map $1_R:R\to R$ factors through a localized polynomial
$K$-algebra, i.e. there is an element $f\in K[x_1,\ldots,x_n]$ the
polynomial ring over $K$ and there are $K$-algebra homomorphisms
$\varphi:R\to K[x_1,\ldots,x_n][1/f]$ and
$\psi:K[x_1,\ldots,x_n][1/f]\to R$ satisfying
$\psi\circ\varphi=1_R$.
\end{defn}

It is not difficult to see that ``rational" $\Rightarrow$ ``stably rational" $\Rightarrow$ ``retract rational".

One of the motivation to study Noether's problem arises from the
inverse Galois problem. If $K$ is an infinite field, it is known
that $K(G)$ is retract rational over $K$ if and only if there
exists a generic Galois $G$-extension over $K$ [\relax {Sa1},
Theorem 5.3; \relax {Sa3}, Theorem 3.12], which guarantees the
existence of a Galois $G$-extension of $K$, provided that $K$ is a
Hilbertian field. On the other hand, the existence of a generic
Galois $G$-extension over $K$ is equivalent to the existence of a
generic polynomial for $G$. For the relationship among these
notions, see \cite{DM}. For a survey of Noether's problem the
reader is referred to articles of Swan and Kersten \cite{Sw2,Ke}.

Although Noether's problem for abelian groups was investigated
extensively, our knowledge for the non-abelian Noether's problem
was amazingly scarce (see, for example, \cite{Ka2}). We will list
below some previous results of non-abelian Noether's problem,
which are relevant to the theme of this article.

\begin{theorem}[Saltman \cite{Sa2}]
For any prime number $p$ and for any field $K$ with char $K\ne p$
$($in particular, $K$ may be an algebraically closed field$)$,
there is a meta-abelian $p$-group $G$ of order $p^9$ such that $K(G)$ is not retract rational over $K$.
In particular, it is not rational.
\end{theorem}

\begin{theorem}[Hajja \cite{Ha}]
Let $G$ be a finite group containing an abelian normal subgroup $N$ such that $G/N$ is a cyclic group of order $<23$.
Then ${\Bbb C}(G)$ is rational over ${\Bbb C}$.
\end{theorem}

\begin{theorem}[Chu and Kang \cite{CK}]
Let $p$ be a prime number, $G$ be a $p$-group of order $\le p^4$
with exponent $p^e$. Let $K$ be any field such that either char
$K=p$ or char $K\ne p$ and $K$ contains a primitive $p^e$-th root
of unity. Then $K(G)$ is rational over $K$.
\end{theorem}

\begin{theorem}[Kang \cite{Ka4}]
Let $G$ be a metacyclic $p$-group with exponent $p^e$,
and let $K$ be any field such that {\rm (i)} char $K=p$, or {\rm (ii)} char $K\ne p$ and $K$ contains a primitive
$p^e$-th root of unity.
Then $K(G)$ is rational over $K$.
\end{theorem}

Note that, in Theorems 1.3---1.5, it is assumed that the ground
field contains enough roots of unity. We may wonder whether ${\Bbb
Q}(G)$ is rational if $G$ is a non-abelian $p$-group of small
order. The answer is rather optimistic when $G$ is a group of
order 8 or 16.

\begin{theorem}[\cite{CHK,Ka3}]
Let $K$ be any field and $G$ be any non-abelian group of order 8 or 16 other than the generalized quaternion group of order 16.
Then $K(G)$ is always rational over $K$.
\end{theorem}

However Serre was able to show that ${\Bbb Q}(G)$ is not rational when $G$ is the generalized quaternion group
[\relax{Se}, p.441--442; \relax{GMS}, Theorem 33.26 and Example 33.27, p.89--90].

On the other hand, if $p$ is an odd prime number, Saltman proves the following theorem.

\begin{theorem}[Saltman \cite{Sa4}]
Let $p$ be an odd prime number and $G$ be a non-abelian group of order $p^3$.
If $K$ is a field containing a primitive $p$-th root of unity,
then $K(G)$ is stably rational.
\end{theorem}

The above theorem may be generalized to the case of $p$-groups containing a maximal cyclic subgroup, namely,

\begin{theorem}[Hu and Kang \cite{HuK}]
Let $p$ be a prime number and $G$ be a non-abelian $p$-group of order $p^n$ containing a cyclic subgroup of index $p$.
If $K$ is any field containing a primitive $p^{n-2}$-th root of unity,
then $K(G)$ is rational over $K$.
\end{theorem}

In this article we will prove the following theorem.

\begin{theorem}
Let $p$ be an odd prime number, $G$ be the non-abelian group of
exponent $p$ and of order $p^3$ or $p^4$. If $K$ is a field with
$[K(\zeta_p):K] \leq 2$, then $K(G)$ is rational over $K$.
\end{theorem}

The rationality problem of $K(G)$ seems rather intricate if the
ground field $K$ has no enough root of unity. We don't know the
answer to the rationality of $K(G)$ when the assumption that
$[K(\zeta_p):K] \leq 2$ is waived in the above theorem. On the
other hand, as to the retract rationality of $K(G)$, a lot of
information may be obtained. Before stating our results, we recall
a theorem of Saltman first.

\begin{theorem}[{Saltman \cite[Theorem 3.5]{Sa1}}]
Let $K$ be a field, $G=A \rtimes G_0$ be a semi-direct product
group where $A$ is an abelian normal subgroup of $G$. Assume that
$gcd \{|A|,|G_0| \}=1$ and both $K(A)$ and $K(G_0)$ are retract
rational over $K$. Then $K(G)$ is retract rational over $K$.
\end{theorem}

Thus the main problem is to investigate the retract rationality
for $p$-groups. We will prove $K(G)$ is retract rational for many
$p$-groups $G$ of exponent $p$.

\begin{theorem}
Let $p$ be a prime number, $K$ be any field, and $G=A \rtimes G_0$
be a semi-direct product group where $A$ is a normal elementary
$p$-group of $G$ and $G_0$ is a cyclic group of order $p^m$. If
$p=2$ and $char K \neq 2$, assume furthermore that
$K(\zeta_{2^m})$ is a cyclic extension of $K$. Then $K(G)$ is
retract rational over $K$.
\end{theorem}

If $p$ is an odd prime number, a $p$-group of exponent $p$
containing an abelian normal subgroup of index $p$ certainly
satisfies the assumption in Theorem 1.11. In particular, a
$p$-group of exponent $p$ and of order $p^3$ or $p^4$ belongs to
this class of $p$-groups (see \cite{CK}).  There are six
$p$-groups of exponent $p$ and of order $p^5$; only four of them
contain abelian normal subgroups of index $p$. Previously the
retract rationality of $K(G)$ for non-abelian $p$-groups, i.e. the
existence of generic polynomials for such groups $G$, is known
only when $G$ is of order $p^3$ and of exponent $p$.

Similarly if $G=A \rtimes G_0$ is a semi-direct product of
$p$-groups such that $A$ is a normal subgroup of order $p$, and
$G_0$ is a direct product of an elementary $p$-group with a cyclic
group of order $p^m$, then $G$ also satisfies the assumption in
Theorem 1.11 provided that the assumption that $K(\zeta_{2^m})$ is
a cyclic extension of $K$ remains in force.

The above Theorem 1.11 is deduced from the following theorem.

\begin{theorem}
Let $K$ be any field, and $G=A \rtimes G_0$ be a semi-direct
product group where $A$ is a normal abelian subgroup of exponent
$e$ and $G_0$ is a cyclic group of order $m$. Assume that

{\rm (i)} either $char K =0$ or $char K > 0$ with $char K \nmid
em$, and

{\rm (ii)} both $K(\zeta_e)$ and $K(\zeta_m)$ are cyclic
extensions of $K$ such that $gcd \{m, [K(\zeta_e):K] \}=1$.

Then $K(G)$ is retract rational over $K$.
\end{theorem}

The idea of the proof of Theorem 1.12 is to add a primitive $e$-th
root of unity to the ground field and the question is reduced to a
question of multiplicative group actions. It is Voskresenskii who
realizes that the multiplicative group action is related to the
birational classification of algebraic tori \cite{Vo}. However,
the multiplicative group action arising in the present situation
is not the function field of an algebraic torus; it is a new type
of multiplicative group actions. Thus we need a new criterion for
retract rationality. It is the following theorem.

\begin{theorem}
Let $\pi_1$ and $\pi_2$ be finite abelian groups, $\pi=\pi_1\times \pi_2$, and $L$
be a Galois extension of the field $K$ with $\pi_1={Gal}({L}/{K})$. Regard $L$ as a
$\pi$-field through the projection $\pi\to \pi_1$. Assume that

{\rm (i)} $\gcd \{|\pi_1|,|\pi_2|\}=1$,

{\rm (ii)} char $K=0$ or $char K >0$ with $char K \nmid |\pi_2|$,
and

{\rm (iii)} $K(\zeta_m)$ is a cyclic extension of $K$ where $m$ is
the exponent of $\pi_2$.

If $M$ is $\pi$-lattice such that $\rho_\pi(M)$ is an invertible
$\pi$-lattice, than $L(M)^\pi$ is retract rational over $K$.
\end{theorem}

The reader will find that the above theorem is an adaptation of
Saltman's criterion for retract rational algebraic tori
\cite[Theorem 3.14]{Sa3} (see Theorem 2.5). We also formulate
another criterion for retract rationality of $L(M)^\pi$ when $\pi$
is a semi-direct product group (see Theorem 4.3). An amusing
consequence of this criterion (when compared with Theorem 1.3) is
that, if $G=A\rtimes H$ is a semi-direct product of an abelian
normal subgroup $A$ and a cyclic subgroup $H$, then ${\Bbb C}(G)$
is always retract rational (see Proposition 5.2).

We will organize this paper as follows. We recall some basic facts
of multiplicative group actions in Section 2. In particular, the
flabby class map which was mentioned in Theorem 1.13 will be
defined. We will give additional tools for proving Theorem 1.9 and
Theorems 1.11--1.13 in Section 3. In Section 4 Theorem 1.13 and
its variants will be proved. The proof of Theorem 1.11 and Theorem
1.12 will be given in Section 5. Section 6 contains the proof of
Theorem 1.9.

Acknowledgements. I am indebted to Prof. R. G. Swan for providing
a simplified proof in Step 7 of Case 1 of Theorem 1.9 (see Section
6). The proof in a previous version of this paper was lengthy and
complicated. I thank Swan's generosity for allowing me to include
his proof in this article.

Notations and terminology. A field extension $L$ over $K$ is
rational if $L$ is purely transcendental over $K$; $L$ is stably
rational over $K$ if there exist $y_1,\ldots,y_N$ such that
$y_1,\ldots,y_N$ are algebraically independent over $L$ and
$L(y_1,\ldots,y_N)$ is rational over $K$. More generally, two
fields $L_1$ and $L_2$ are called stably isomorphic if
$L_1(x_1,\ldots,x_m)$ is isomorphic to $L_2(y_1,\ldots,y_n)$ where
$x_1,\ldots, x_m$ and $y_1,\ldots,y_n$ are algebraically
independent over $L_1$ and $L_2$ respectively.

Recall the definition of $K(G)$ at the beginning of this section:
$K(G)=K(x_g:g\in G)^G$. If $L$ is a field with a finite group $G$
acting on it, we will call it a $G$-field. Two $G$-fields $L_1$
and $L_2$ are $G$-isomorphic if there is an isomorphism
$\varphi:L_1\to L_2$ satisfying $\varphi(\sigma\cdot
u)=\sigma\cdot\varphi(u)$ for any $\sigma\in G$, any $u\in L_1$.

We will denote by $\zeta_n$ a primitive $n$-th root of unity in
some extension field of $K$ when char $K=0$ or char $K=p>0$ with
$p\nmid n$. All the groups in this article are finite groups.
${\Bbb Z}_n$ will be the cyclic group of order $n$ or the ring of
integers modulo $n$ depending on the situation from the context.
${\Bbb Z}[\pi]$ is the group ring of a finite group $\pi$ over
${\Bbb Z}$. $Z(G)$ is the center of the group $G$. The exponent of
a group $G$ is the least common multiple of the orders of elements
in $G$. The representation space of the regular representation of
$G$ over $K$ is denoted by $W=\bigoplus_{g\in G}K\cdot x(g)$ where
$G$ acts on $W$ by $g\cdot x(h)=x(gh)$ for any $g,h\in G$.

\section{Multiplicative group actions}

Let $\pi$ be a finite group.
A $\pi$-lattice $M$ is a finitely generated ${\Bbb Z}[\pi]$-module such that $M$ is a free abelian group when it is
regarded as an abelian group.

For any field $K$ and a $\pi$-lattice $M$,
$K[M]$ will denote the Laurent polynomial ring and $K(M)$ is the quotient field of $K[M]$.
Explicitly, if $M=\bigoplus_{1\le i\le m} {\Bbb Z}\cdot x_i$ as a free abelian group,
then $K[M]=K[x^{\pm 1}_1,\ldots,x^{\pm 1}_m]$ and $K(M)=K(x_1,\ldots,x_m)$.
Since $\pi$ acts on $M$,
it will act on $K[M]$ and $K(M)$ by $K$-automorphisms,
i.e.\ if $\sigma\in \pi$ and $\sigma\cdot x_j=\sum_{1\le i\le m} a_{ij} x_i\in M$,
then we define the action of $\sigma$ in $K[M]$ and $K(M)$ by $\sigma\cdot x_j=\prod_{1\le i\le m} x^{a_{ij}}_i$.

The multiplicative action of $\pi$ on $K(M)$ is called a purely
monomial action in \cite{HK1}. If $\pi$ is a group acting on the
rational function field $K(x_1,\ldots,x_m)$ by $K$-automorphism
such that $\sigma\cdot x_j=c_j(\sigma)\cdot\prod_{1\le i\le m}
x^{a_{ij}}_i$ where $\sigma\in\pi$, $a_{ij}\in{\Bbb Z}$ and $c_j
(\sigma)\in K\backslash \{0\}$, such a multiplicative group action
is called a monomial action. Monomial actions arise when studying
Noether's problem for non-split extension groups \cite{Ha,Sa5}.

We will introduce another kind of multiplicative actions. Let
$K\subset L$ be fields and $\pi$ be a finite group. Suppose that
$\pi$ acts on $L$ by $K$-automorphisms (but it is not assumed that
$\pi$ acts faithfully on $L$). Given a $\pi$-lattice $M$, the
action of $\pi$ on $L$ can be extended to an action of $\pi$ on
$L(M)$ ($=L(x_1,\ldots,x_m)$ if $M=\bigoplus_{1\le i\le m}{\Bbb
Z}\cdot x_i$) by $K$-automorphisms defined as follows: If
$\sigma\in\pi$ and $\sigma\cdot x_j=\sum_{1\le i\le m}a_{ij}x_i\in
M$, then the multiplication action in $L(M)$ is defined by
$\sigma\cdot x_j=\prod_{1\le i\le m}x^{a_{ij}}_i$ for $1\le j\le
m$.

When $L$ is a Galois extension of $K$ and $\pi={Gal}({L}/{K})$
(and therefore $\pi$ acts faithfully on $L$), the fixed subfield
$L(M)^\pi$ is the function field of the algebraic torus defined
over $K$, split by $L$ and with character group $M$ (see
\cite{Vo}).

We recall some basic facts of the theory of flabby (flasque)
$\pi$-lattices developed by Endo and Miyata, Voskresenskii,
Colliot-Th\'{e}l\`{e}ne and Sansuc, etc. \cite{Vo,CTS}. We refer
the reader to \cite{Sw2,Sw3,Lo} for a quick review of the theory.

In the sequel, $\pi$ denotes a finite group unless otherwise specified.

\begin{defn}\rm
A $\pi$-lattice $M$ is called a permutation lattice if $M$ has a
${\Bbb Z}$-basis permuted by $\pi$. $M$ is called an invertible
(or permutation projective)  lattice, if it is a direct summand of
some permutation lattice. A $\pi$-lattice $M$ is called a flabby
(or flasque) lattice if $H^{-1}(\pi',M)=0$ for any subgroup $\pi'$
of $\pi$. (Note that $H^{-1}(\pi',M)$ denotes the Tate cohomology
group.) Similarly, $M$ is called coflabby if $H^1(\pi',M)=0$ for
any subgroup $\pi'$ of $\pi$.
\end{defn}

It is known that an invertible $\pi$-lattice is necessarily a
flabby lattice [\relax{Sw2}, Lemma 8.4; \relax{Lo}, Lemma 2.5.1].

\begin{theorem}[Endo and Miyata {[\relax{Sw3}, Theorem 3.4; \relax{Lo}, 2.10.1]}]
Let $\pi$ be a finite group. Then any flabby $\pi$-lattice is
invertible if and only if all Sylow subgroups of $\pi$ are cyclic.
\end{theorem}

Denote by ${\cal L}_\pi$ the class of all $\pi$-lattices,
and by ${\cal F}_\pi$ the class of all flabby $\pi$-lattices.

\begin{defn}\rm
We define an equivalence relation $\sim$ on ${\cal F}_\pi$:
Two $\pi$-lattices $E_1$ and $E_2$ are similar,
denoted by $E_1\sim E_2$,
if $E_1\oplus P$ is isomorphic to $E_2\oplus Q$ for some permutation lattices $P$ and $Q$.
The similarity class containing $E$ will be denoted by $[E]$.
Define $F_\pi={\cal F}_\pi/\sim$,
the set of all similarity classes of flabby $\pi$-lattices.
\end{defn}

$F_\pi$ becomes a commutative monoid if we define $[E_1]+[E_2]=[E_1\oplus E_2]$.
The monoid $F_\pi$ is called the flabby class monoid of $\pi$.

\begin{defn}\rm
We define a map $\rho:{\cal L}_\pi\to F_\pi$ as follows. For any
$\pi$-lattice $M$, there exists a flabby resolution, i.e. a short
exact sequence of $\pi$-lattices $0\to M\to P\to E \to 0$ where
$P$ is a permutation lattice and $E$ is a flabby lattice
\cite[Lemma 8.5]{Sw2}. We define $\rho_\pi(M)=[E]\in F_\pi$. The
map $\rho_\pi:{\cal L}_\pi\to F_\pi$ is well-defined \cite[Lemma
8.7]{Sw2}; it is called the flabby class map. We will simply write
$\rho$ instead of $\rho_\pi$, if the group $\pi$ is obvious from
the context.
\end{defn}

\begin{theorem}[{Saltman \cite[Theorem 3.14]{Sa3}}]
Let $L$ be a Galois extension of $K$ with $\pi={Gal}({L}/{K})$ and $M$ be a
$\pi$-lattice. Then $\rho_\pi(M)$ is invertible if and only if $L(M)^\pi$ is retract
rational over $K$.
\end{theorem}

\section{Generalities}

We recall several results which will be used later.

\begin{theorem}[{\cite[Theorem 1]{HK2}}]
Let $L$ be a field and $G$ be a finite group acting on $L(x_1,\ldots,x_m)$,
the rational function field of $m$ variables over $L$.
Suppose that

{\rm (i)} for any $\sigma\in G$, $\sigma(L)\subset L$;

{\rm (ii)} The restriction of the action of $G$ to $L$ is faithful;

{\rm (iii)} for any $\sigma\in G$,
\[
\left(\begin{array}{c} \sigma(x_1) \\ \vdots \\ \sigma(x_m) \end{array}\right)=
A(\sigma) \left(\begin{array}{c} x_1 \\ \vdots \\ x_m \end{array}\right) +B(\sigma)
\]
where $A(\sigma)\in GL_m(L)$ and $B(\sigma)$ is an $m\times 1$ matrix over $L$.

Then $L(x_1,\ldots,x_m)=L(z_1,\ldots,z_m)$ where $\sigma(z_i)=z_i$
for all $\sigma\in G$, and for any $1\le i \le m$. In fact,
$z_1,\ldots,z_m$ can be defined by
\[
\left(\begin{array}{c} z_1 \\ \vdots \\ z_m \end{array}\right)=
A\cdot \left(\begin{array}{c} x_1 \\ \vdots \\ x_m \end{array}\right) +B
\]
for some $A\in GL_m(L)$ and for some $B$ which is an $m\times 1$
matrix over $L$. Moreover, if $B(\sigma)=0$ for all $\sigma\in G$,
we may choose $B=0$ in defining $z_1,\ldots,z_m$.
\end{theorem}

\begin{theorem}[{Kuniyoshi \cite[Theorem 2.5]{CHK}}]
Let $K$ be a field with $char K=p>0$ and $G$ be a $p$-group. Then
$K(G)$ is always rational over $K$.
\end{theorem}

\begin{prop}
Let $\pi$ be a finite group and $L$ be a $\pi$-field.
Suppose that $0\to M_1\to M_2\to N\to 0$ is a short exact sequence of $\pi$-lattices satisfying
{\rm (i)} $\pi$ acts faithfully on $L(M_1)$,
and {\rm (ii)} $N$ is an invertible $\pi$-lattice.
Then the $\pi$-fields $L(M_2)$ and $L(M_1\oplus N)$ are $\pi$-isomorphic.
\end{prop}

{\it Proof}. We follow the proof of {[\relax {Le}, Proposition
1.5]}. Denote $L(M_1)^\times=L(M_1)\backslash \{0\}$. Consider the
exact sequence of $\pi$-modules:
\[
0\to L(M_1)^\times\to L(M_1)^\times\cdot M_2\to N\to 0.
\]

By Hilbert Theorem 90, we find that $H^1(\pi', L(M_1)^\times)=0$
for any subgroup $\pi' \subset \pi$. Applying {[\relax {Le},
Proposition 1.2]} we find that the above exact sequence splits.
The resulting $\pi$-morphism $N\to L(M_1)^\times \cdot M_2$
provides the required $\pi$-isomorphism form $L(M_2)$ to
$L(M_1\oplus N)$.\hfill$\square$
\begin{lemma}
Let the assumptions be the same as in Proposition 3.3. Assume
furthermore that $N$ is a permutation $\pi$-lattice. Then
$L(M_2)^{\pi}$ is rational over $L(M_1)^{\pi}$.
\end{lemma}

{\it Proof}.
By Proposition 3.3, $L(M_2)=L(M_1)(N)$.
Since $\pi$ acts faithfully on $L(M_1)$,
we may apply Theorem 3.1 and find $u_1,\ldots,u_n\in L(M_2)$ such that $L(M_2)=L(M_1)(u_1$, $\ldots,u_n)$ with
$\sigma(u_i)=u_i$ for any $\sigma\in \pi$,
any $1\le i\le n$ where $n=\fn{rank}(N)$.
Hence $L(M_2)^\pi=L(M_1)^\pi(u_1,\ldots,u_n)$.\hfill$\square$

\begin{lemma}
Let $\pi$ be a finite abelian group of exponent $e$ and $K$ be a
field such that char $K=0$ or $char K > 0$ with $char K \nmid e$.
If $P$ is a permutation $\pi$-lattice, then
$K(P)^\pi=K(\zeta_e)(M)^{\pi_0}$ where
$\pi_0={Gal}({K(\zeta_e)}/{K})$ and $M$ is some $\pi_0$-lattice.
\end{lemma}

{\it Proof}. We follow the standard approach to solving Noether's
problem for abelian groups {[\relax {Sw1; Sw2; Le}]}.

Note that $K(P)^\pi=\{K(\zeta_e)(P)^{\langle
\pi_0\rangle}\}^{\langle\pi\rangle}=K(\zeta_e)(P)^{\langle\pi,\pi_0\rangle}$
where the action of $\pi$ is extended to $K(\zeta_e)(P)$ by
defining $g(\zeta_e)=\zeta_e$ for any $g\in \pi$, and the action
of $\pi_0$ is extended to $K(\zeta_e)(P)$ by requiring that
$\pi_0$ acts trivially on $P$.

Since $\pi$ is abelian of exponent $e$, we may diagonalize its
action on $P$, i.e. we may find $x_1,\ldots,x_n\in K(\zeta_e)(P)$
such that $n=\fn{rank}(P)$, $g{(x_i)}/{x_i}\in
\langle\zeta_e\rangle$ for any $g \in \pi$, and
$K(\zeta_e)(P)=K(\zeta_e)(x_1,\ldots,x_n)$.

Thus
$K(\zeta_e)(P)^{\langle\pi\rangle}=K(\zeta_e)(y_1,\ldots,y_n)$
where $y_1,\ldots,y_n$ are monomials in $x_1,\ldots,x_n$. Let $M$
be the multiplicative subgroup generated by $y_1,\ldots,y_n$ in
$K(\zeta_e)(y_1,\ldots,y_n)\backslash\{0\}$. Then $M$ is a
$\pi_0$-lattice and
$K(\zeta_e)(y_1,\ldots,y_n)^{\pi_0}=K(\zeta_e)(M)^{\pi_0}$.\hfill$\square$

\begin{prop}
Let $\pi$ and $K$ be the same as in Lemma 3.5. If $K(\zeta_e)$ is
a cyclic extension of $K$, then $K(\pi)$ is retract rational over
$K$.
\end{prop}

{\it Proof}. We may regard the regular representation of $\pi$ is
given by a permutation $\pi$-lattice. Thus
$K(\pi)=K(\zeta_e)(M)^{\pi_0}$ where
$\pi_0={Gal}({K(\zeta_e)}/{K})$. Since $\pi_0$ is assumed cyclic,
thus we may apply Theorem 2.2 and Theorem 2.5.\hfill$\square$

\section{Proof of Theorem 1.13}

\begin{lemma}
Let $\pi$ be a finite group, $M$ be a $\pi$-lattice. Suppose that $\pi_0$ is a
normal subgroup of $\pi$ and $\pi_0$ acts trivially on $M$. Thus we may regard $M$
as a lattice over ${\pi}/{\pi_0}$.

{\rm (1)} $M$ is a permutation $\pi$-lattice $\Leftrightarrow$ So is it as a
${\pi}/{\pi_0}$-lattice.

{\rm (2)} $M$ is an invertible $\pi$-lattice $\Leftrightarrow$ So is it as a
${\pi}/{\pi_0}$-lattice.

{\rm (3)} $M$ is a flabby $\pi$-lattice $\Leftrightarrow$ So is it as a
${\pi}/{\pi_0}$-lattice.

{\rm (4)} If $0\to M\to P\to E\to 0$ is a flabby resolution of $M$ as a
${\pi}/{\pi_0}$-lattice, this short exact sequence is also a flabby resolution of
$M$ as a $\pi$-lattice.

{\rm (5)} $\rho_\pi(M)$ is an invertible $\pi$-lattice $\Leftrightarrow
\rho_{{\pi}/{\pi_0}}(M)$ is an invertible ${\pi}/{\pi_0}$-lattice.
\end{lemma}

{\it Proof}. The properties (1)--(4) can be found in \cite[Lemma
2, p.179--180]{CTS}. As to (5), the direction ``$\Leftarrow$" is
obvious by applying (4). For the other direction, assume
$\rho_\pi(M)$ is an invertible $\pi$-lattice. Let $0\to M\to P\to
E\to 0$ be a flabby resolution of $M$ as a $\pi$-lattice. Then
$0\to M^{\pi_0}\to P^{\pi_0}\to E^{\pi_0}\to 0$ is a flabby
resolution of $M=M^{\pi_0}$ in the category of
${\pi}{/\pi_0}$-lattices by \cite[Lemma (xi), p.180]{CTS}. It
remains to show that $E^{\pi_0}$ is invertible. Since
$[E]=\rho_\pi(M)$ is invertible, we can find a $\pi$-lattice $N$
such that $E\oplus N=N'$ is a permutation $\pi$-lattice. Note that
${N'}^{\pi_0}$ is a permutation ${\pi}/{\pi_0}$-lattice by
\cite[Lemma 2(i), p.180]{CTS}. We find that $E^{\pi_0}$ is
invertible because $E^{\pi_0}\oplus N^{\pi_0}=(E\oplus
N)^{\pi_0}={N'}^{\pi_0}$. \hfill$\square$

\begin{lemma}
Let the assumptions be the same as in Theorem 1.13. If $P$ is a
permutation $\pi$-lattice, then $L(P)^\pi$ is retract rational
over $K$.
\end{lemma}

{\it Proof}. Since $\pi_1$ and $\pi_2$ are abelian groups with
$\gcd\{|\pi_1|,|\pi_2|\}=1$, every subgroup $\pi'$ of $\pi$ can be
written as $\pi'=\rho\times \lambda$ where $\rho$ is a subgroup of
$\pi_1$ and $\lambda$ is a subgroup of $\pi_2$.

As a permutation $\pi$-lattice, we may write $P=\bigoplus {\Bbb
Z}[{\pi}/{\pi^{(i)}}]$ where $\pi^{(i)}$ is a subgroup of $\pi$.
Write $\pi^{(i)}=\rho_i\times \lambda_i$ where $\rho_i \subset
\pi_1$, $\lambda_i \subset \pi_2$. Hence ${\Bbb Z}[\pi /
\pi^{(i)}]= {\Bbb Z}[\pi/(\rho_i \times \lambda_i)]={\Bbb
Z}[(\pi_1/\rho_i) \times (\pi_2/\lambda_i)]$.

It is not difficult to see that ${\Bbb
Z}[{\pi}/{\pi^{(i)}}]=\bigoplus{\Bbb Z}\cdot u^{(i)}_{kl}$ where
$1\le k\le t=|{\pi_1}/{\rho_i}|$, $1\le l\le r=
|{\pi_2}/{\lambda_i}|$. Moreover, if $g\in \pi_1$ and
$g'\in\pi_2$, then $g\cdot u^{(i)}_{kl}=u^{(i)}_{g(k),l}$,
$g'\cdot u^{(i)}_{kl}=u^{(i)}_{k,g'(l)}$ and the homomorphisms
$\pi_1\to S_t$ and $\pi_2\to S_r$ are induced from the permutation
representations associated to ${\pi}/{\pi^{(i)}}=(\pi_1/\rho_i)
\times (\pi_2/\lambda_i)$ where $S_t$ and $S_r$ are the symmetric
groups of degree $t$ and $r$ respectively.

Since $\pi_1$ is faithful on $L$, we may apply Theorem 3.1.
Explicitly, for any $1\le l\le r$, we may find $A^{(i)}\in
GL_t(L)$ and define $v^{(i)}_{kl}$ by
\begin{eqnarray}
\left(\begin{array}{c}  v^{(i)}_{1,l} \\ v^{(i)}_{2,l} \\ \vdots \\ v^{(i)}_{t,l}  \end{array}\right) =
A^{(i)}\left(\begin{array}{c}  u^{(1)}_{1,l} \\ u^{(i)}_{2,l} \\ \vdots \\ u^{(i)}_{t,l}  \end{array}\right)
\label{eq1}
\end{eqnarray}
such that $g\cdot v^{(i)}_{k,l}=v^{(i)}_{k,l}$ for any $g\in
\pi_1$, and $L({\Bbb Z}[{\pi}/{\pi^{(i)}}])=L(u^{(i)}_{kl}: 1\le
k\le t,~1\le l\le r$. If $g'\in \pi_2$, from the relation $g'\cdot
u^{(i)}_{kl}=u^{(i)}_{k,g'(l)}$ and Formula (\ref{eq1}), we find
that $g'\cdot v^{(i)}_{kl}=v^{(i)}_{k,g'(l)}$.

Since
$L(P)^{\pi_1}=L(u^{(i)}_{kl})^{\pi_1}=L(v^{(i)}_{kl})^{\pi_1}=K(v^{(i)}_{kl})$
where $i$, $k$, $l$ runs over index sets which are understood, it
follows that $L(P)^\pi=K(v^{(i)}_{kl})^{\pi_2}$. Note that $\pi_2$
acts on $\{v^{(i)}_{kl}\}$ by permutations.

By Lemma 3.5 $K(v^{(i)}_{kl})^{\pi_2}=K(\zeta_m)(M)^{\pi_0}$ where
$m$ is the exponent of $\pi_2$, $\pi_0={Gal}({K(\zeta_m)}/{K})$
and $M$ is some $\pi_0$-lattice.

By our assumption, $\pi_0$ is a cyclic group. Hence
$\rho_{\pi_0}(M)$ is invertible by Theorem 2.2. Apply Theorem 2.5.
We find that $K(\zeta_m)(M)^{\pi_0}$ is retract rational over $K$.
\hfill$\square$\\

{\it Proof of Theorem} 1.13 -----------------------------------

Step 1. Suppose that $M$ is a $\pi$-lattice such that
$\rho_\pi(M)$ is invertible.

Define $\pi_0=\{g\in\pi_2:g$ acts trivially on $L(M)\}$. Then
${\pi}/{\pi_0}$ acts faithfully on $L(M)$. Moreover,
$\rho_{{\pi}/{\pi_0}}(M)$ is invertible by Lemma 4.1.

In other words, without loss of generality we may assume that
$\pi$ is faithfully on $L(M)$ . Thus we will keep in force this
assumption in the sequel.

Step 2. Since $\rho_\pi(M)$ is invertible, by \cite[Theorem 2.3, p.176]{Sa3},
we may find $\pi$-lattice $M'$, $P$, $Q$ such that $P$ and $Q$ are permutation lattices,
$0\to M\to M'\to Q\to 0$ is exact,
and the inclusion map $M\to M'$ factors through $P$,
i.e. the following diagram commutes
\begin{eqnarray}
\begin{array}{l}
\put(-25,17){$M$} \put(-12,13){\vector(1,-1){10}}
\put(-39,-6){$1_M$} \put(-21,11){\vector(0,-1){30}}  \put(1,-8){$P$} \put(11,-9){\vector(1,-1){10}}
\put(-55,-33){0}\put(-45,-29){\vector(1,0){15}} \put(-25,-33){$M$} \put(-11,-29){\vector(1,0){28}}\put(22,-33){$M'$}
\put(41,-29){\vector(1,0){15}} \put(61,-33){$Q$} \put(74,-29){\vector(1,0){15}}\put(95,-33){0.}
\end{array} \label{eq2}
\end{eqnarray}

The remaining proof proceeds quite similar to that of
\cite[Theorem 3.14, p.189]{Sa3}.

Step 3. We get a commutative diagram of $K$-algebra morphisms from the diagram in (\ref{eq2}), i.e.
\begin{eqnarray*}
\put(-70,32){$L[M]^\pi$} \put(-42,24){\vector(2,-1){20}}
\put(-67,3){$1$} \put(-56,24){\vector(0,-1){38}} \put(-21,1){$L[P]^\pi$} \put(5,-5){\vector(2,-1){20}}
\put(-70,-30){$L[M]^\pi$} \put(-34,-27){\vector(1,0){48}} \put(21,-30){$L[M']^\pi$}
\put(59,-27){\vector(1,0){20}} \put(85,-30){$L[Q]^\pi$}
\end{eqnarray*}

Step 4. The quotient field of $L[P]^\pi$ is $L(P)^\pi$, which is
retract rational over $K$ by Lemma 4.2. Thus the identity map
$1:L[P]^\pi\to L[P]^\pi$ factors rationally by \cite[Lemma
3.5]{Sa3}, i.e. there is a localized polynomial ring
$K[x_1,\ldots,x_n][1/f]$ and $K$-algebra maps $\varphi:L[P]^\pi\to
K[x_1,\ldots,x_n][1/f]$, $\psi:K[x_1,\ldots,x_n][1/f]\to L[P]^\pi$
such that $\psi\circ\varphi =1$.

It follows that the composite map $g: L[M]^\pi\to L[P]^\pi\to
L[M']^\pi$ also factors rationally, i.e. there are $K$-algebra
$\varphi':L[M]^\pi\to K[x_1,\ldots,x_n][1/f]$,
$\psi':K[x_1,\ldots,x_n][1/f]\to L[M']^\pi$ such that
$g=\psi'\circ \varphi'$.

Step 5. By Lemma 3.4 $L(M')^\pi$ is rational over $L(M)^\pi$.
(This is the only one step we use the assumption that $\pi$ is faithful on $L(M)$.)

Now we may apply \cite[Proposition 3.6(b), p.183]{Sa3} where, in
the notation of \cite{Sa3}, we take $S=T=L[M]^\pi$, $\varphi$ is
the identity map on $L(M)^\pi$. We conclude that $1:L[M]^\pi\to
L[M]^\pi$ factors rationally,
i.e. $L(M)^\pi$ is retract rational over $K$.\hfill$\square$\\

Here is a variant of Theorem 1.13.

\begin{theorem}
Let $\pi$ be a finite group, $0\to \pi_1\to \pi \to \pi_2 \to 1$
is a group extension, and $L$ be a Galois extension of the field
$K$ with $\pi_2={Gal}({L}/{K})$. Let $\pi$ act on $L$ through the
projection $\pi\to \pi_2$. Assume that

{\rm (i)} $\pi_1$ is an abelian group of exponent $e$ with
$\zeta_e \in L$,

{\rm (ii)} the extension $0\to\pi_1\to\pi\to\pi_2\to 1$ splits,
and

{\rm (iii)} every Sylow subgroup of $\pi_2$ is cyclic.

If $M$ is a $\pi$-lattice such that $\rho_\pi(M)$ is an invertible
lattice, then $L(M)^\pi$ is retract rational over $K$.
\end{theorem}

{\it Proof}. The proof is very similar to the proof of Theorem
1.13.

We claim that $L(P)^\pi$ is retract rational for any permutation
$\pi$-lattice $P$.

For the proof, we will use \cite[Theorem 2.1, p.546]{Sa5}. We will
show that (c) of \cite[Theorem 2.1]{Sa5} is valid, which will
guarantee that $L(P)^\pi$ is retract rational. By assumption
(iii), part (d) of \cite[Theorem 2.1]{Sa5} is valid by Theorem
2.2. It remains to check that the embedding problem of ${L}/{K}$
and the extension $0\to\pi_1\to\pi\to\pi_2\to 1$ is solvable. But
this is the well-known split embedding problem \cite[Theorem 1.9,
p.12]{ILF}.

Now define $\pi_0=\{g\in\pi:g$ acts trivially on $L(M)\}$. Note
that $\pi_0 \subset \pi_1$. The remaining proof is the same as in
the proof of Theorem 1.13 and is omitted.\hfill$\square$

\begin{coro}
Let $\pi$ be an abelian group of exponent $e$, $K$ be a field with
$\zeta_e\in K$. Suppose that $M$ is a $\pi$-lattice and $\pi$ acts
on $K(M)$ by $K$-automorphisms. If $\rho_\pi(M)$ is an invertible
module, than $K(M)^\pi$ is retract rational over $K$.
\end{coro}

\section{Proof of Theorems 1.11 and 1.12}

{\it Proof of Theorem 1.12} ----------------------------------

\

Step 1.  Write $G= A \rtimes G_0$ where $G_0=<\sigma>$ is a cyclic
group of order $m$, and $A=A_1 \times A_2 \times \cdots \times
A_r$ with each $A_i=<\rho_i>$ being a cyclic group of order $e_i$,
$e= e_1$ and $e_r \ | \ e_{r-1}, \ \cdots, e_2 \ | \ e_1$.

Define $B_i= < \rho_i, \cdots, \hat {\rho_i}, \cdots, \rho_r>$ to be the subgroup of
$A$ by deleting the $i$-th direct factor.

Let $Gal(K(\zeta_e)/K)=<\tau>$ be a cyclic group of order $n$.
Write $\zeta=\zeta_e$ and $\tau \cdot \zeta = \zeta^a$ for some
integer $a$.

For $1 \le i \le r$, choose $\zeta_i \in <\zeta>$ such that
$\zeta_i$ is a primitive $e_i$-th root of unity.  Note that $\tau
\cdot \zeta_i=\zeta_i^a$.

Let $n_i=[K(\zeta_i):K]$ for $1 \le i \le r$.  Note that $n=n_1$.

\

Step 2.  Let $W=\bigoplus_{g \in G} K \cdot x(g)$ be the
representation space of the regular representation of $G$.  Then
$$K(G)=K(W)^G= \{K(\zeta)(W)^{<\tau>} \}^G= K(\zeta)(W)^{<G, \tau>}$$
where the action of $G$ and $\tau$ are extended to $K(\zeta)(W)$
by requiring that $G$ and $\tau$ act trivially on $K(\zeta)$ and
$W$ respectively.

\

Step 3. For $1 \le i \le r$, define $$u_i= \sum_{l\in {\Bbb
Z}_{e_i}, \ g \in B_i}\zeta_i^l \cdot x(\rho_i^l \cdot g).$$

Note that $\rho_i \cdot u_i = \zeta_i^{-1}u_i$ and, if $j \ne i$,
then $\rho_i \cdot u_i=u_i$.

For $1 \le i \le r$, write $\chi_i$ to be the character $\chi_i: A
\longrightarrow K(\zeta)^\times$ such that $g \cdot u_i= \chi_i(g)
\cdot u_i$.  Since $\chi_1, \cdots, \chi_r$ are distinct
characters on $A$, it follows that $u_1, \cdots, u_r$ are linearly
independent vectors in $K(\zeta) \otimes_KW$.

Moreover, the subgroup $A$ acts faithfully on $\bigoplus _{1 \le i
\le r} K(\zeta)u_i$.

\

Step 4.  For $1 \le i \le r$, $j \in {\Bbb Z}_{n_i}$, define
$$v_{ij}= \sum_{l\in {\Bbb Z}_{e_i}, \ g \in B_i}\zeta_i^{la^j}
\cdot x(\rho_i^l \cdot g).$$

Thus we have elements $v_{1,1}, \cdots, v_{1,n_1}, v_{2,1},
\cdots, v_{2,n_2}, \cdots, v_{r,n_r}$; in total we have $n_1+n_2+
\cdots +n_r$ such elements.

As in Step 3, these elements $v_{ij}$ are linearly independent.  Note that $\tau
\cdot v_{ij}=v_{i, j+1}$.

\

Step 5.  For $1 \le i \le r$, $j \in {\Bbb Z}_{n_i}, \ k \in {\Bbb
Z}_m$.  define

$$x_{i,j,k}= \sum_{l\in {\Bbb Z}_{e_i}, \ g \in B_i}\zeta_i^{la^j } \cdot x(\sigma^k
 \cdot \rho_i^l  \cdot g).$$

We have in total $m(n_1+n_2+ \cdots +n_r)$ such elements
$x_{i,j,k}$.

Note that $x_{i, j, 0}=v_{ij}$, and the vector $x_{i, j, k}$ is
just a ``translation'' of the vector $x_{i, j, 0}$ in the space
$\bigoplus_{g \in G} K(\zeta)x(g)$ (with basis $x(g), \ g \in G$).
Thus these vectors $x_{i, j, k}$ are linearly independent.  Note
that $\sigma \cdot v_{i, j, k}=v_{i, j, k+1}$. Moreover, the group
$<G, \tau>$ acts faithfully on $\bigoplus K(\xi)\cdot x_{i, j,
k}$.

Apply Theorem 3.1.  We find $K(G)=K(\zeta)(x_{ijk}: \ 1 \le i \le
r, \ j \in {\Bbb Z}_{n_i}, k \in {\Bbb Z}_m)(w_1, \cdots, w_s)$
where $s=|G|-m(n_1+n_2+ \cdots +n_r)$, and $\lambda(w_i)=w_i$ for
any $i$, any $\lambda \in <G, \tau>$.

\

Step 6.  We will consider the fixed field $K(\zeta)(x_{ijk}: 1 \le
i \le r, \ j \in {\Bbb Z}_{n_i}, k \in {\Bbb Z}_m)^{<G, \tau>}$.

Let $\pi=<\sigma, \tau>$, $\Lambda= {\Bbb Z}[\pi]$.  Let
$N=<x_{ijk}:1 \le i \le r, \ j \in {\Bbb Z}_{n_i}, k \in {\Bbb
Z}_m>$ be the multiplicative subgroup generated by these $x_{ijk}$
in $K(\zeta)(x_{ijk}: 1 \le i \le r, \ j \in {\Bbb Z}_{n_i}, k \in
{\Bbb Z}_m)\setminus \{0\}$.  Note that $N$ is a $\Lambda$-module.
In fact, $N$ is a permutation $\pi$-lattice.

Define $$\Phi: N \longrightarrow {\Bbb Z}_{n_1} \times {\Bbb
Z}_{n_2} \times \cdots \times {\Bbb Z}_{n_r}$$ by $\Phi(x)=(\bar
l_1, \bar l_2, \cdots, \bar l_r) \in {\Bbb Z}_{n_1} \times {\Bbb
Z}_{n_2} \times \cdots \times {\Bbb Z}_{n_r}$, if $x=\prod
x_{ijk}^{b_{ijk}}$ (with $b_{ijk} \in {\Bbb Z}$) and $\rho_1(x)=
\zeta_1^{l_1}x$, $\rho_2(x)=\zeta_2^{l_2}x$, $\cdots$,
$\rho_r(x)=\zeta_r^{l_r}x$.

Define $M=\fn{Ker} (\Phi)$.  We find that $K(\zeta)(x_{ijk}: 1 \le
i \le r, j \in {\Bbb Z}_{n_i}, k \in {\Bbb Z}_m)^A= K(\zeta)(M)$.
Thus it remains to find $K(\zeta)(M)^{\pi}$.  Note that $M$ is a
$\pi$-lattice.

\

Step 7.  Since $gcd \{m, \ n\}=1$, it follows that $\pi$ is a
cyclic group.  Hence $\rho_{\pi}(M)$ is invertible by Theorem 2.2.

Apply Theorem 1.13.  We find that $K(\zeta)(M)^{\pi}$ is retract
rational over $K$. Since $K(G)$ is rational over
$K(\zeta)(M)^{\pi}$, it follows that $K(G)$ is also retract
rational over $K$ by \cite[Proposition 3.6(a), p.183]{Sa3}. \hfill
$\square$

\

{\it Proof of Theorem 1.11} -------------------------------

If $char K = p$, apply Theorem 3.2. Thus $K(G)$ is rational. In
particular it is retract rational.

>From now on we will assume that $char K \ne p$. It is not
difficult to verify that all the assumptions of Theorem 1.12 are
valid in the present situation. Hence the result. \hfill $\square$

\begin{coro}  Let $K$ be a field, $p$ be any prime number, $G=A\rtimes G_0$ be a
semi-direct product of $p$-groups where $A$ is a normal abelian
subgroup of exponent $p^e$ and $G_0$ is a cyclic group of order
$p^m$.  If $char K \ne p$, assume that {\rm (i)} both
$K(\zeta_{p^e})$ and $K(\zeta_{p^m})$ are cyclic extensions of $K$
and {\rm (ii)} $p \nmid [K(\zeta_{p^e}):K]$.  Then $K(G)$ is
retract rational over $K$.
\end{coro}

\begin{prop} Let $K$ be a field and $G=A \rtimes G_0$ be a semi-direct product.
Assume that

{\rm (i)} $A$ is an abelian normal subgroup of exponent $e$ in $G$ and $G_0$ is a
cyclic group of order $m$.

{\rm (ii)} either $char K=0$ or $char K >0$ with $char K \nmid
em$, and

{\rm (iii)} $\zeta_e \in K$ and $K(\zeta_m)$ is a cyclic extension
of $K$.

If $G \longrightarrow GL(V)$ is a finite-dimensional linear representation of $G$
over $K$, then $K(V)^G$ is retract rational over $K$. \end{prop}

{\it Proof}.  Decompose $V$ into $V=\bigoplus V_{\chi}$ where
$\chi: A \longrightarrow K\setminus \{0\}$ runs over linear
characters of $A$ and $V_{\chi}=\{v \in V: \ g \cdot
v=\chi(g)\cdot v \ \fn{ for \ all} \ g \in A\}$

It is easy to see that $\sigma(V_{\chi})=V_{\chi^{\prime}}$ for
some $\chi^{\prime}$. Suppose $V_{\chi_1}, V_{\chi_2}, \cdots,
V_{\chi_s}$ is an orbit of $\sigma$, i.e. $\sigma(V_{\chi_j}) =
V_{\chi_{j+1}}$ and $\sigma(V_{\chi_s})= V_{\chi_1}$. Choose a
basis $v_1, \cdots, v_t$ of $V_{\chi_1}$.  It follows that
$\{\sigma^j(v_i): \ 1 \le i \le t, \ 0 \le j \le s-1\}$ is a basis
of $\oplus _{1 \le j \le s}V_{\chi_j}$.

In this way, we may find a basis $w_1, \cdots, w_n$ in $V$ such
that (i) for any $g \in A$, any $1 \le i \le r$,  \, $g \cdot
w_i=\alpha w_i$ for some $\alpha \in K$; and (ii) $\sigma \cdot
w_i=w_j$ for some $j$.

It follows that $K(V)^A= K(w_1, \cdots, w_n)^A=K(u_1, \cdots,
u_n)$ where $u_1, \cdots, u_n$ are monomials in $w_1, \cdots,
w_n$.  Thus $K(u_1, \cdots, u_n)=K(M)$ for some $G_0$-lattice $M$.
It follows that $K(V)^G=\{K(V)^A\}^{G_0}=K(M)^{<\sigma>}$.

By Theorem 2.2, \, $\rho_{G_0}(M)$ is invertible.  Apply Theorem
1.12. \hfill $\square$

\begin{remark} \rm In the above theorem it is essential that $G$ is a
is a split extension group.  For non-split extension groups,
monomial actions (instead of merely purely monomial actions) may
arise; see the proof of Theorem 1.3 \cite{Ha} and also
\cite{Sa5}.\end{remark}

\section{Proof of Theorem 1.9}

We will prove Theorem 1.9 in this section.

If $char K = p$, apply Theorem 3.2. We find that $K(G)$ is
rational. In particular, it is retract rational.

Thus we will assume that $char K \neq p$ from now on till the end
of this section.

If $\zeta_p \in K$, we apply Theorem 1.4 to find that $K(G)$ is
rational. Hence we will consider only the situation $[K(\zeta_p):
K]=2$ in the sequel.

Since $p$ be an odd prime number, there is only one non-abelian
$p$-group of order $p^3$ with exponent $p$, and there are
precisely two non-isomorphic non-abelian $p$-groups of order $p^4$
with exponent $p$ (see \cite[Section 2 and Section 3]{CK}). We
will solve the rationality problem of these three groups
separately.

Let $\zeta=\zeta_p$, $\fn{Gal}(K(\zeta)/K)=\langle \tau\rangle$
with $\tau\cdot \zeta=\zeta^{-1}$.

\begin{case}\rm
$G=\langle \sigma_1,\sigma_2,\sigma_3:\sigma^p_i=1,~\sigma_1\in
Z(G),~\sigma_2\sigma_3=\sigma_3\sigma_1\sigma_2 \rangle$.
\end{case}

Step 1. Let $W=\bigoplus_{g\in G} K\cdot x(g)$ be the
representation space of the regular representation of $G$. Note
that
\[
K(G)=K(x(g):g\in G)^G=\{K(\zeta)(x(g):g\in G)^{\langle
\tau\rangle}\}^G=K(\zeta)(x(g):g\in G)^{\langle G,\tau\rangle}.
\]

\

Step 2. For $i\in{\Bbb Z}_p$, define
$x_{0,i},x_{1,i}\in\bigoplus_{g\in G} K(\zeta)\cdot x(g)$ by
\begin{eqnarray*}
x_{0,i}=\sum_{j,k\in{\Bbb Z}_p} \zeta^{-j-k} x(\sigma^i_3\sigma^j_1\sigma^k_2), \\
x_{1,i}=\sum_{j,k\in{\Bbb Z}_p} \zeta^{j+k}\cdot
x(\sigma^i_3\sigma^j_1\sigma^k_2).
\end{eqnarray*}

We find that
\begin{eqnarray*}
\sigma_1: &\hspace*{-4mm}& x_{0,i}\mapsto \zeta x_{0,i},\ x_{1,i}\mapsto \zeta^{-1} x_{1,i}, \\
\sigma_2: &\hspace*{-4mm}& x_{0,i}\mapsto \zeta^{i+1} x_{0,i},\ x_{1,i}\mapsto \zeta^{-i-1} x_{1,i}, \\
\sigma_3: &\hspace*{-4mm}& x_{0,i}\mapsto x_{0,i+1},\ x_{1,i}\mapsto x_{1,i+1}, \\
\tau:     &\hspace*{-4mm}& x_{0,i}\leftrightarrow x_{1,i}.
\end{eqnarray*}

The restriction of the action of $\langle\sigma_1,\sigma_2
\rangle$ to $K(\zeta)\cdot x_{0,i}$ and $K(\zeta)\cdot x_{1,i}$ is
given by distinct characters of $\langle\sigma_1,\sigma_2 \rangle$
to $K(\zeta) \setminus \{ 0 \}$. Hence
$x_{0,0},x_{0,1},x_{0,2},\ldots,x_{0,p-1}, x_{1,0},\ldots,
x_{1,p-1}$ are linearly independent vectors. Moreover, the action
of $\langle G,\tau \rangle$ on $K(\zeta)(x_{0,i},x_{1,i}:i\in{\Bbb
Z}_p)$ is faithful.

By Theorem 3.1 $K(\zeta)(x(g):g\in
G)=K(\zeta)(x_{0,i},x_{1,i}:i\in{\Bbb Z}_p)(X_1,\ldots,X_l)$ where
$l=p^3-2p$ and $\rho(X_j)=X_j$ for any $1\le j \le l$, any
$\rho\in \langle G,\tau\rangle$.

\

Step 3. Define $x_0=x^p_{0,0}$, $y_0=x_{0,0} x_{1,0}$ and
$x_i=x_{0,i}\cdot x_{0,i-1}^{-1}$, $y_i=y_{0,i}\cdot
y^{-1}_{0,i-1}$ for $1\le i\le p-1$. It follows that
$K(\zeta)(x_{0,i},x_{1,i}:i\in{\Bbb Z}_p)^{\langle\sigma_1
\rangle}=K(\zeta)(x_i,y_i:i\in{\Bbb Z}_p)$. Moreover, the actions
of $\sigma_2,\sigma_3,\tau$ are given as
\begin{eqnarray*}
\sigma_2: &\hspace*{-4mm}& x_0\mapsto x_0,\ y_0\mapsto y_0,\
x_i\mapsto \zeta x_i,\ y_i\mapsto \zeta^{-1}y_i
                           \mbox{ ~for~ }1\le i\le p-1, \\
\sigma_3: &\hspace*{-4mm}& x_0\mapsto x_0x^p_1,\ x_1\mapsto x_2\mapsto \cdots\mapsto x_{p-1}\mapsto(x_1x_2\cdots x_{p-1})^{-1}, \\
          &\hspace*{-4mm}& y_0\mapsto y_0y_1x_1,\ y_1\mapsto y_2\mapsto \cdots \mapsto y_{p-1}\mapsto(y_1y_2\cdots y_{p-1})^{-1}, \\
\tau: &\hspace*{-4mm}& x_0\mapsto y^p_0x^{-1}_0,\ y_0\mapsto y_0,\
x_i\leftrightarrow y_i \mbox{ ~for ~} 1\le i\le p-1.
\end{eqnarray*}

\

Step 4. Define $X=x_0y^{-(p-1)/2}_0$, $Y=x^{-1}_0 y^{(p+1)/2}_0$.

Then $K(\zeta)(x_i,y_i: i\in{\Bbb Z}_p)=K(\zeta)(X,Y,x_i,y_i:1\le
i\le p-1)$ and $\sigma_2(X)=X$, $\sigma_2(Y)=Y$,
$\sigma_3(X)=\alpha X$, $\sigma_3(Y)=\beta Y$, $\tau:X
\leftrightarrow Y$ where $\alpha,\beta\in K(\zeta)(x_i,y_i:1\le
i\le p-1)\backslash \{0\}$.

Apply Theorem 3.1. We may find $\widetilde{X}$, $\widetilde{Y}$ so
that $K(\zeta)(x_i,y_i:i\in{\Bbb Z}_p)= K(\zeta)(x_i,y_i:1\le i\le
p-1)(\widetilde{X},\widetilde{Y})$ with
$\rho(\widetilde{X})=\widetilde{X}$,
$\rho(\widetilde{Y})=\widetilde{Y}$ for any $\rho\in \langle
\sigma_2,\sigma_3,\tau\rangle$. Thus it remains to consider
$K(\zeta)(x_i,y_i:1\le i\le p-1)^{\langle
\sigma_2,\sigma_3,\tau\rangle}$.

\

Step 5. Define $u_0=x^p_1$, $v_0=x_1y_1$, $u_i=x_{i+1}x^{-1}_i$,
$v_i=y_{i+1}y^{-1}_i$ for $1\le i\le p-2$.

It follows that $K(\zeta)(x_i,y_i:1\le i\le p-1)^{\langle
\sigma_2\rangle}=K(\zeta)(u_i,v_i:0\le i\le p-2)$. The actions of
$\sigma_3$ and $\tau$ are given by
\begin{eqnarray*}
\sigma_3: &\hspace*{-4mm}& u_0\mapsto u_0u^p_1,\ v_0\mapsto
v_0v_1u_1,\
              u_1\mapsto u_2\mapsto\cdots\mapsto u_{p-2}\mapsto (u_0u^{p-1}_1u^{p-2}_2\cdots u^2_{p-2})^{-1}, \\
          &\hspace*{-4mm}& v_1\mapsto v_2\mapsto\cdots\mapsto v_{p-2}\mapsto u_0(v^p_0v^{p-1}_1\cdots v^2_{p-2})^{-1}, \\
\tau: &\hspace*{-4mm}& u_0\mapsto u^{-1}_0v^p_0,\ v_0\mapsto v_0,\
u_i\leftrightarrow v_i \mbox{ ~for ~} 1\le i\le p-2.
\end{eqnarray*}

\

Step 6. Define $u_{p-1}=(u_0u^{p-1}_1u^{p-2}_2\cdots
u^2_{p-2})^{-1}$, $w_i=v_0v_1\cdots v_i u_1u_2\cdots u_i$ for
$1\le i\le p-2$, and $w_{p-1}=(v^{p-1}_0v^{p-2}_1\cdots v_{p-2}
u^{p-2}_1u^{p-3}_2\cdots u_{p-2})^{-1}$. We find that
$K(u_i,v_i:0\le i\le p-2)=K(u_i,w_i:1\le i\le p-1)$ and
\begin{eqnarray*}
\sigma_3: &\hspace*{-4mm}& u_1\mapsto u_2\mapsto\cdots\mapsto u_{p-1}\mapsto (u_1u_2\cdots u_{p-1})^{-1}, \\
          &\hspace*{-4mm}& w_1\mapsto w_2\mapsto\cdots\mapsto w_{p-1}\mapsto (w_1w_2\cdots w_{p-1})^{-1}, \\
\tau: &\hspace*{-4mm}& u_i\mapsto w_i(u_iw_{i-1})^{-1},\
w_i\mapsto w_i \mbox{ ~for ~} 1\le i\le p-1.
\end{eqnarray*}
where we write $w_0=v_0$ for convenience. (Granting that $w_0$ is
defined as above, we have a relation $w_0w_1\cdots w_{p-1}=1$. But
we don't have the relation $u_0u_1\cdots u_{p-1}=1$ because we
define $u_0,u_1,\ldots,u_{p-2}$ first and $u_{p-1}$ is defined by
another way.)

We will study whether $K(\zeta)(u_i,w_i:1\le i\le p-1)^{\langle
\sigma_3,\tau\rangle}$ is rational or not.

\

Step 7. The multiplicative action in Step 6 can be formulated in
terms of $\pi$-lattices where $\pi=\langle \tau,\sigma_3\rangle$
as follows.

Let $M=(\bigoplus_{1\le i\le p-1} {\Bbb Z}\cdot
u_i)\oplus(\bigoplus_{1\le i\le p-1}{\Bbb Z}\cdot w_i)$ and define
$w_0=-w_1-w_2-\cdots-w_{p-1}$. Define a ${\Bbb Z}[\pi]$-module
structure on $M$ by
\begin{eqnarray*}
\sigma_3: &\hspace*{-4mm}& u_1\mapsto u_2\mapsto\cdots\mapsto u_{p-1}\mapsto -u_1-u_2-\cdots-u_{p-1}, \\
          &\hspace*{-4mm}& w_1\mapsto w_2\mapsto\cdots\mapsto w_{p-1}\mapsto -w_1-w_2-\cdots-w_{p-1}, \\
\tau: &\hspace*{-4mm}& u_i\mapsto -u_i+w_i-w_{i-1},\ w_i\mapsto
w_i \mbox{ ~for ~} 1\le i\le p-1.
\end{eqnarray*} .

We claim that $M\simeq {\Bbb Z}[\pi_1]\bigotimes_{\Bbb Z} {\Bbb
Z}[\pi_2]/\Phi_p(\sigma_3)$ with $\pi_1=\langle\tau\rangle$,
$\pi_2=\langle\sigma_3\rangle$.

Throughout this step, we will write $\sigma = \sigma_3$ and
$\pi=<\sigma, \ \tau>$. Define $\rho = \sigma \tau$. Then $\rho$
is a generator of $\pi$ with order $2p$ where $p$ is an odd prime
number.

Let $\Phi_p(T) \in {\Bbb Z}[T]$ be the $p$-th cyclotomic
polynomial.  Note that
$\Phi_p(\sigma)(u_i)=\Phi_p(\sigma)(w_i)=0$. It follows that
$\Phi_p(\sigma) \cdot M=0$.

Since $\Phi_p(\sigma^2)=\sum_{0 \le i \le p-1}\sigma^i=
\Phi_p(\sigma)$, we know that $\Phi_p(\sigma^2)\cdot M=0$.  From
$\rho^2=\sigma^2$, we find that $\Phi_p(\rho^2) \cdot M=0$.  It is
well-known that $\Phi_p(T^2)=\Phi_p(T) \Phi_{2p}(T)$ (see [Ka1,
Theorem 1.1] for example).  We conclude that
$\Phi_p(\rho)\Phi_{2p}(\rho)\cdot M=0$. In other words, we may
regard $M$ as a module over $\Lambda={\Bbb
Z}[\pi]/<\Phi_p(\rho)\Phi_{2p}(\rho)>$.

Clearly $\Lambda$ is isomorphic to ${\Bbb
Z}[\pi_1]\bigotimes_{\Bbb Z} {\Bbb Z}[\pi_2]/\Phi_p(\sigma)$ where
$\pi_1=\langle\tau\rangle$, $\pi_2=\langle\sigma\rangle$.

It remains to show that $M$ is isomorphic to $\Lambda$ as a
$\Lambda$-module.

It is not difficult to verify that
$\Phi_{2p}(\rho)(u_1)=f(\rho)(w_1)$ where $f(T)=\Phi_{2p}(T) -
T^{p-1}$. Define $v=u_1-w_1$. We find that
$\Phi_{2p}(\rho)(v)=-w_0$.

Since $M$ is a $\Lambda$-module generated by $u_1$ and $w_0$, it
follows that, as a $\Lambda$-module, $M=< u_1, w_0> = <u_1- \sigma
w_0, w_0>$ $=<u_1-w_1, w_0>=<v, w_0>$ $=<v,
w_0+\Phi_{2p}(\rho)(v)>=<v,w_0-w_0>=<v>$, i.e. $M$ is a cyclic
$\Lambda$-module generated by $v$. Thus we get an epimorphism
$\Lambda \rightarrow M$, which is an isomorphism by counting the
$\Bbb Z$-ranks of both sides. Hence the result.

\

Step 8. By Step 7, $M\simeq {\Bbb Z}[\pi_1]\bigotimes_{\Bbb Z}
{\Bbb Z}[\pi_2]/\Phi_p(\sigma_3)$ where
$\pi_1=\langle\tau\rangle$, $\pi_2=\langle\sigma_3\rangle$.

Thus we may choose a ${\Bbb Z}$-basis $Y_1,Y_2,\ldots,Y_{p-1}$,
$Z_1,\ldots,Z_{p-1}$ for $M$ such that
\begin{eqnarray*}
\sigma_3: &\hspace*{-4mm}& Y_1\mapsto Y_2\mapsto\cdots\mapsto Y_{p-1}\mapsto -Y_1-\cdots-Y_{p-1}, \\
          &\hspace*{-4mm}& Z_1\mapsto Z_2\mapsto\cdots\mapsto Z_{p-1}\mapsto -Z_1-\cdots-Z_{p-1}, \\
\tau: &\hspace*{-4mm}& Y_i\leftrightarrow Z_i \mbox{ ~for ~} 1\le
i\le p-1.
\end{eqnarray*}

Hence $K(\zeta)(M)=K(\zeta)(Y_i,Z_i:1\le i\le p-1)$. We emphasize
that $\sigma_3$ acts on $Y_i$, $Z_i$ by multiplicative actions on
the field $K(\zeta)(M)$ and $\sigma_3\cdot \zeta=\zeta$,
$\tau\cdot \zeta=\zeta^{-1}$.

\

Step 9. In the field $K(\zeta)(M)$, define
\[
\begin{array}{l}
s_0=1+Y_1+Y_1Y_2+Y_1Y_2Y_3+\cdots+Y_1Y_2\cdots Y_{p-1},\\
t_0=1+Z_1+Z_1Z_2+\cdots+Z_1Z_2\cdots Z_{p-1}, \\
s_1=(1/s_0)-(1/p),\ s_2=(Y_1/s_0)-(1/p),\ \ldots,\ s_{p-1}=(Y_1Y_2\cdots Y_{p-2}/s_0)-(1/p), \\
t_1=(1/t_0)-(1/p),\ t_1=(Z_1/t_0)-(1/p),\ \ldots,\
t_{p-1}=(Z_1Z_2\cdots Z_{p-2}/t_0)-(1/p).
\end{array}
\]

It is easy to verify that $K(\zeta)(M)=K(\zeta)(s_i,t_i:1\le i\le
p-1)$ and
\begin{eqnarray*}
\sigma_3: &\hspace*{-4mm}& s_1\mapsto s_2\mapsto\cdots\mapsto s_{p-1}\mapsto -s_1-s_2-\cdots-s_{p-1}, \\
          &\hspace*{-4mm}& t_1\mapsto t_2\mapsto\cdots\mapsto t_{p-1}\mapsto -t_1-t_2-\cdots-t_{p-1}, \\
\tau: &\hspace*{-4mm}& s_i \leftrightarrow t_i.
\end{eqnarray*}

\

Step 10. Define $r_i=s_i+t_i$ for $1\le i\le p-1$. Then
$K(\zeta)(s_i,t_i:1\le i\le p-1)=K(\zeta)(s_i,r_i:1\le i\le p-1)$.
Note that
\begin{eqnarray*}
\sigma_3: &\hspace*{-4mm}& r_1\mapsto r_2\mapsto\cdots\mapsto r_{p-1}\mapsto -r_1-r_2-\cdots-r_{p-1}, \\
\tau: &\hspace*{-4mm}& r_i\mapsto r_i,\ s_i\mapsto -s_i+r_i.
\end{eqnarray*}

Apply Theorem 3.1. We find that $K(\zeta)(s_i,r_i:1\le i\le
p-1)=K(\zeta)(r_i:1\le i\le p-1)(A_1,\ldots,A_{p-1})$ for some
$A_1,\ldots,A_{p-1}$ with $\sigma_3(A_i)=\tau(A_i)=A_i$ for $1\le
i\le p-1$.

Thus $K(\zeta)(r_i,s_i:1\le i\le
p-1)^{\langle\tau\rangle}=K(\zeta)(r_i:1\le i\le
p-1)^{\langle\tau\rangle}(A_1,\ldots,A_{p-1})
=K(r_1,\ldots,r_{p-1},A_1,\ldots,A_{p-1})$.

It remains to find $K(r_1,\ldots,r_{p-1})^{\langle
\sigma_3\rangle}$.

\

Step 11. Write $\pi_2=\langle\sigma_3\rangle$. The $\pi_2$-fields
$K(r_1,\ldots,r_{p-1},A_1)$ and $K(B_0,B_1,\ldots,B_{p-1})$ are
$\pi_2$-isomorphic where $\sigma_3:B_0\mapsto
B_1\mapsto\cdots\mapsto B_{p-1}\mapsto B_0$. For, we may define
$B=B_0+B_1+\cdots+B_{p-1}$ and $C_i=B_i-(B/p)$ for $0\le i\le
p-1$.

In other words
$K(r_1,\ldots,r_{p-1},A_1,\ldots,A_{p-1})^{\langle\sigma_3\rangle}=
K(B_0,B_1,\ldots,B_{p-1})^{\langle\sigma_3\rangle}(A_2,\ldots$,
$A_{p-1})=K({\Bbb Z}_p)(A_2,\ldots,A_{p-1})$.

By Lemma 3.5, $K({\Bbb Z}_p)=K(\zeta)(N)^{\pi_1}$ where
$\zeta=\zeta_p$, $\pi_1=\fn{Gal}(K(\zeta)/K)=\langle\tau\rangle$
with $\tau\cdot \zeta=\zeta^{-1}$, and $N$ is some
$\pi_1$-lattice.

By Reiner's Theorem \cite{Re}, the $\pi_1$-lattice $N$ is a direct
sum of lattices of three types: (1) $\tau:z\mapsto -z$, (2)
$\tau:z\mapsto z$, (3) $\tau:z_1\leftrightarrow z_2$. Thus we find
$K(\zeta)(N)^{\langle \tau\rangle} =K(\zeta)(z_1$,
$\ldots,z_a,z'_1,z'_2,\ldots,z'_b,
z''_1,w''_1,\ldots,z''_c,w''_c)$ where $\tau:z_1\mapsto 1/z_1$,
$\ldots$, $z_a\mapsto 1/z_a$, $z'_1\mapsto z'_1$, $\ldots$,
$z'_b\mapsto z'_b$, $z''_1\leftrightarrow w''_1$, $\ldots$, $z''_c
\leftrightarrow z''_c$.

By Theorem 3.1 we may ``neglect" the roles of
$z'_1,\ldots,z'_b,z''_1,w''_1,\ldots,z''_c,w''_c$. We may
linearize the action on $z_1,\ldots,z_a$ by defining
$w_i=1/(1+z_i)$ when $1\le i\le a$. Then $\tau:w_i\mapsto -w_i+1$.
Thus we may ``neglect" the roles of $w_i$ by Theorem 3.1 again. In
conclusion,
$K(\zeta)(z_1,\ldots,z_a,z'_1,\ldots,z'_b,z''_1,w''_1,\ldots,z''_c,w''_c)^{\langle\tau\rangle}$
is rational over $K$.\hfill$\square$

\begin{case} \rm
$G=\langle \sigma_1,\sigma_2,\sigma_3,\sigma_4:\sigma^p_i=1,\
\sigma_1,\sigma_2\in Z(G),\
\sigma^{-1}_4\sigma_3\sigma_4=\sigma_1\sigma_3\rangle$.
\end{case}

The proof is very similar to Case 1.

Step 1. For $i\in{\Bbb Z}_p$, define $x_{0,i}, x_{1,i}, y_{0,i},
y_{1,i}\in \bigoplus_{g\in G} K(\zeta)\cdot x(g)$ by
\begin{eqnarray*}
x_{0,i}&=&\sum_{j,k,l\in{\Bbb Z}_p} \zeta^{-j-k-l} x(\sigma^i_4\sigma^j_3\sigma^k_1\sigma^l_2), \\
x_{1,i}&=&\sum_{j,k,l\in{\Bbb Z}_p} \zeta^{j+k+l} x(\sigma^i_4\sigma^j_3\sigma^k_1\sigma^l_2), \\
y_{0,i}&=&\sum_{j,k,l\in{\Bbb Z}_p} \zeta^{-j-k+l} x(\sigma^i_4\sigma^j_3\sigma^k_1\sigma^l_2), \\
y_{1,i}&=&\sum_{j,k,l\in{\Bbb Z}_p} \zeta^{j+k-l}
x(\sigma^i_4\sigma^j_3\sigma^k_1\sigma^l_2).
\end{eqnarray*}

The action of $\langle G,\tau\rangle$ is given by
\begin{eqnarray*}
\sigma_1: &\hspace*{-4mm}& x_{0,i}\mapsto\zeta x_{0,i},\
x_{1,i}\mapsto\zeta^{-1}x_{1,i},\
                           y_{0,i}\mapsto\zeta y_{0,i},\ y_{1,i}\mapsto\zeta^{-1}y_{1,i}, \\
\sigma_2: &\hspace*{-4mm}& x_{0,i}\mapsto\zeta x_{0,i},\
x_{1,i}\mapsto\zeta^{-1}x_{1,i},\
                           y_{0,i}\mapsto\zeta^{-1}y_{0,i},\ y_{1,i}\mapsto\zeta y_{1,i}, \\
\sigma_3: &\hspace*{-4mm}& x_{0,i}\mapsto\zeta^{i+1}x_{0,i},\
x_{1,i}\mapsto\zeta^{-i-1}x_{1,i},\
                           y_{0,i}\mapsto\zeta^{i+1}y_{0,i},\ y_{1,i}\mapsto\zeta^{-i-1}y_{1,i}, \\
\sigma_4: &\hspace*{-4mm}& x_{0,i}\mapsto x_{0,i+1},\
x_{1,i}\mapsto x_{1,i+1},\
                           y_{0,i}\mapsto y_{0,i+1},\ y_{1,i}\mapsto y_{1,i+1}, \\
\tau: &\hspace*{-4mm}& x_{0,i}\leftrightarrow x_{1,i},\
y_{0,i}\leftrightarrow y_{1,i}.
\end{eqnarray*}

Note that $x_{0,i}$, $x_{1,i}$, $y_{0,i}$, $y_{1,i}$ where
$i\in{\Bbb Z}_p$ are linearly independent vectors in
$\bigoplus_{g\in G}K(\zeta)\cdot x(g)$. Apply Theorem 3.1. It
suffices to consider the rationality problem of
$K(\zeta)(x_{0,i},x_{1,i},y_{0,i},y_{1,i}: i\in{\Bbb
Z}_p)^{\langle G,\tau\rangle}$.

\

Step 2. Define $x_0=x_{0,0}^p$, $y_0=x_{0,0}x_{1,0}$,
$X_0=y_{0,0}(x_{0,0})^{-1}$, $Y_0=y_{1,0}(x_{1,0})^{-1}$; for $1
\leq i \leq p-1$, define $x_i=x_{0,i}(x_{0,i-1})^{-1}$,
$y_i=x_{1,i}(x_{1,i-1})^{-1}$, $X_i=y_{0,i}(y_{0,i-1})^{-1}$,
$Y_i=y_{1,i}(y_{1,i-1})^{-1}$. Then
$K(\zeta)(x_{0,i},x_{1,i},y_{0,i},y_{1,i}:i\in{\Bbb
Z}_p)^{\langle\sigma_1\rangle}= K(\zeta)(x_i$,
$y_i,X_i,Y_i:i\in{\Bbb Z}_p)$ and the actions are given by
\begin{eqnarray*}
\sigma_2: &\hspace*{-4mm}& x_0\mapsto x_0,\ y_0\mapsto y_0,\ X_0\mapsto\zeta^{-2}X_0,\ Y_0\mapsto\zeta^2 Y_0, \\
          &\hspace*{-4mm}& \mbox{All the other generators are fixed by ~}\sigma_2; \\
\sigma_3: &\hspace*{-4mm}& x_0\mapsto x_0,\ y_0\mapsto y_0,\ X_0\mapsto X_0,\ Y_0\mapsto Y_0, \\
          &\hspace*{-4mm}& x_i\mapsto\zeta x_i,\ y_i\mapsto\zeta^{-1}y_i,\ X_i\mapsto\zeta X_i,\ Y_i\mapsto\zeta^{-1}Y_i, \\
\sigma_4: &\hspace*{-4mm}& x_0\mapsto x_0x^p_1,\ y_0\mapsto
y_0y_1x_1,\
                           X_0\mapsto X_0X_1x^{-1}_1,\ Y_0\mapsto Y_0Y_1y^{-1}_1, \\
          &\hspace*{-4mm}& x_1\mapsto x_2\mapsto x_3\mapsto\cdots\mapsto x_{p-1}\mapsto(x_1x_2\cdots x_{p-1})^{-1}, \\
          &\hspace*{-4mm}& \mbox{Similarly for ~} y_1,y_2,\ldots,y_{p-1},X_1,\ldots,X_{p-1},Y_1,\ldots,Y_{p-1}; \\
\tau: &\hspace*{-4mm}& x_0\mapsto y^p_0x^{-1}_0,\ y_0\mapsto y_0,\ X_0\mapsto Y_0\mapsto X_0, \\
      &\hspace*{-4mm}& x_i\leftrightarrow y_i,\ X_i\leftrightarrow Y_i \mbox{ ~for~ } 1\le i\le p-1.
\end{eqnarray*}

\

Step 3. Define $\tilde{x}=x_0y^{-(p-1)/2}_0$, $\tilde{y}=x^{-1}_0
y^{(p+1)/2}_0$.

Since $\tilde{x}$, $\tilde{y}$ are fixed by both $\sigma_2$ and
$\sigma_3$, while $\sigma_4:\tilde{x}\mapsto \alpha\tilde{x}$,
$\tilde{y}\mapsto \beta\tilde{y}$,
$\tau:\tilde{x}\leftrightarrow\tilde{y}$ where $\alpha,\beta\in
K(x_1,y_1,X_1,Y_1)\backslash \{0\}$, we may apply Theorem 3.1.
Thus the roles of $\tilde{x}$, $\tilde{y}$ may be ``neglected". It
suffices to consider whether
$K(\zeta)(X_0,Y_0,x_i,y_i,X_i,Y_i:1\le i\le p-1)^
{\langle\sigma_2,\sigma_3,\sigma_4,\tau\rangle}$ is rational over
$K$.

\

Step 4. Define $\widetilde{X}=X^p_0$, $\widetilde{Y}=X_0Y_0$. Then
$K(\zeta)(X_0,Y_0,x_i,y_i,X_i,Y_i:1\le i\le
p-1)^{\langle\sigma_2\rangle}=
K(\zeta)(\widetilde{X},\widetilde{Y},x_i,y_i,X_i,Y_i:1\le i\le
p-1)$. Moreover, the actions on $\widetilde{X}$ and
$\widetilde{Y}$ are given by
\begin{eqnarray*}
\sigma_3: &\hspace*{-4mm}& \widetilde{X}\mapsto \widetilde{X},\ \widetilde{Y}\mapsto \widetilde{Y}, \\
\sigma_4: &\hspace*{-4mm}& \widetilde{X}\mapsto
\widetilde{X}X^p_1x^{-p}_1,\
                           \widetilde{Y}\mapsto \widetilde{Y}X_1Y_1x^{-1}_1y^{-1}_1, \\
\tau: &\hspace*{-4mm}& \widetilde{X}\mapsto
\widetilde{X}^{-1}\widetilde{Y}^p,\
\widetilde{Y}\mapsto\widetilde{Y}.
\end{eqnarray*}

Define $X'=\widetilde{X}\widetilde{Y}^{-(p-1)/2}$,
$Y'=\widetilde{X}^{-1}\widetilde{Y}^{(p+1)/2}$. As in Step 3, we
may apply Theorem 3.1 and ``neglect" the roles of $X'$ and $Y'$.
It remains to make sure whether $K(\zeta)(x_i,y_i,X_i,Y_i:1\le
i\le p-1)^{\langle\sigma_3,\sigma_4,\tau\rangle}$ is rational.

Define $u_0=x^p_1$, $v_0=x_1y_1$; for $1\le i\le p-2$, define
$u_i=x_{i+1}(x^{-1}_i)$, $v_i=y_{i+1}\cdot(y^{-1}_i)$; and for
$1\le i\le p-1$, define $U_i=X_ix^{-1}_i$, $V_i=Y_i y^{-1}_i$. It
follows that $K(\zeta)(x_i,y_i,X_i,Y_i:1\le i\le
p-1)^{\langle\sigma_3\rangle}=
K(\zeta)(u_0,v_0,U_{p-1},V_{p-1},u_i,v_i,U_i,V_i:1\le i\le p-2)$.
Note that the actions of $\sigma_4$ and $\tau$ are given by
\begin{eqnarray*}
\sigma_4: &\hspace*{-4mm}& u_0\mapsto u_0u^p_1,\ v_0\mapsto
v_0v_1u_1,\
                           u_1\mapsto u_2\mapsto\cdots\mapsto u_{p-2}\mapsto(u_0u^{p-1}_1u^{p-2}_2\cdots u^2_{p-2})^{-1}, \\
          &\hspace*{-4mm}& v_1\mapsto v_2\mapsto\cdots\mapsto v_{p-2}\mapsto u_0(v^p_0v^{p-1}_1\cdots v^2_{p-2})^{-1}, \\
          &\hspace*{-4mm}& U_1\mapsto U_2\mapsto\cdots\mapsto U_{p-1}\mapsto(U_1U_2\cdots U_{p-1})^{-1}, \\
          &\hspace*{-4mm}& V_1\mapsto V_2\mapsto\cdots\mapsto V_{p-1}\mapsto(V_1V_2\cdots V_{p-1})^{-1}. \\
\tau: &\hspace*{-4mm}& u_0\mapsto u^{-1}_0v^p_0,\ v_0\mapsto v_0,\ u_i\leftrightarrow v_i\mbox{ ~for~ } 1\le i\le p-2, \\
      &\hspace*{-4mm}& U_i\leftrightarrow V_i \mbox{ ~for ~} 1\le i\le p-1.
\end{eqnarray*}

Note that the actions of $\sigma_4$ and $\tau$ on
$U_1,U_2,\ldots,U_{p-1}$, $V_1,\ldots,V_{p-1}$ may be linearized
by the same method as in Step 9 of Case 1. Hence we may apply
Theorem 3.1 and ``neglect" the roles of $U_i$, $V_i$ for $1\le
i\le p-1$.

In conclusion, all we need is to prove that $K(\zeta)(u_i,v_i:0\le
i\le p-2)^{\langle\sigma_4,\tau\rangle}$ is rational over $K$.

Compare the present situation with that of Step 5 of Case 1. We
have the same generators and the same actions (and even the same
notation). Thus we finish the proof.\hfill$\square$

\begin{case}\rm
$G=\langle \sigma_1,\sigma_2,\sigma_3,\sigma_4:\sigma^p_i=1,\
\sigma_1\in Z(G),\ \sigma_2\sigma_3=\sigma_3\sigma_2,\
\sigma^{-1}_4\sigma_2\sigma_4=\sigma_1\sigma_2,\
\sigma^{-1}_4\sigma_3\sigma_4=\sigma_2\sigma_3\rangle$.
\end{case}

Step 1. For $i\in{\Bbb Z}_p$, define
$x_{0,i},x_{1,i},y_{0,i},y_{1,i}\in \bigoplus_{g\in
G}K(\zeta)\cdot x(g)$ by
\begin{eqnarray*}
x_{0,i}&=&\sum_{j,k,l\in{\Bbb Z}_p} \zeta^{-j-k-l} x(\sigma^i_4\sigma^j_3\sigma^k_2\sigma^l_1), \\
x_{1,i}&=&\sum_{j,k,l\in{\Bbb Z}_p} \zeta^{j+k+l} x(\sigma^i_4\sigma^j_3\sigma^k_2\sigma^l_1), \\
y_{0,i}&=&\sum_{j,k,l\in{\Bbb Z}_p} \zeta^{-j-k+l} x(\sigma^i_4\sigma^j_3\sigma^k_2\sigma^l_1), \\
y_{1,i}&=&\sum_{j,k,l\in{\Bbb Z}_p} \zeta^{j+k-l}
x(\sigma^i_4\sigma^j_3\sigma^k_2\sigma^l_1).
\end{eqnarray*}

The action of $\langle G,\tau\rangle$ is given by
\begin{eqnarray*}
\sigma_1: &\hspace*{-4mm}& x_{0,i}\mapsto\zeta x_{0,i},\
x_{1,i}\mapsto\zeta^{-1}x_{1,i},\
                           y_{0,i}\mapsto\zeta^{-1}y_{0,i},\ y_{1,i}\mapsto\zeta y_{1,i}, \\
\sigma_2: &\hspace*{-4mm}& x_{0,i}\mapsto\zeta^{i+1}x_{0,i},\
x_{1,i}\mapsto\zeta^{-i-1}x_{1,i},\
                           y_{0,i}\mapsto\zeta^{-i+1}y_{0,i},\ y_{1,i}\mapsto\zeta^{i-1}y_{1,i}, \\
\sigma_3: &\hspace*{-4mm}& x_{0,i}\mapsto\zeta^{i+1}x_{0,i},\
x_{1,i}\mapsto\zeta^{-i-1}x_{1,i},\
                           y_{0,i}\mapsto\zeta^{i+1}y_{0,i},\ y_{1,i}\mapsto\zeta^{-i-1}y_{1,i}, \\
\sigma_4: &\hspace*{-4mm}& x_{0,i}\mapsto x_{0,i+1},\
x_{1,i}\mapsto x_{1,i+1},\
                           y_{0,i}\mapsto y_{0,i+1},\ y_{1,i}\mapsto y_{1,i+1}, \\
\tau: &\hspace*{-4mm}& x_{0,i}\leftrightarrow x_{1,i},\
y_{0,i}\leftrightarrow y_{1,i}.
\end{eqnarray*}

Checking the restriction of
$\langle\sigma_1,\sigma_2,\sigma_3\rangle$ as in Step 1 of Case 2,
we find that these vectors $x_{0,i}$, $x_{1,i}$, $y_{0,i}$,
$y_{1,i}$ are linearly independent except possibly for the case
$x_{0,p-2}$ and $y_{1,0}$, and the case $x_{1,p-2}$ and $y_{0,0}$.
But it is easy to see that these vectors are linearly independent,
because their "supports" are disjoint. Apply Theorem 3.1. It
suffices to consider
$K(\zeta)(x_{0,i},x_{1,i},y_{0,i},y_{1,i}:i\in{\Bbb Z}_p)^{\langle
G,\tau\rangle}$.

\

Step 2. Define $x_0=x^p_{0,0}$, $y_0=x_{0,0}x_{1,0}$,
$X_0=x_{0,0}\cdot y_{0,0}$, $Y_0=x_{1,0}y_{1,0}$; and for $1\le
i\le p-1$, define $x_i=x_{0,i}(x_{0,i-1})^{-1}$,
$y_i=x_{1,i}(x_{1,i-1})^{-1}$, $X_i=y_{0,i}(y_{0,i-1})^{-1}$,
$Y_i=y_{1,i}(y_{1,i-1})^{-1}$. Then
$K(\zeta)(x_{0,i},x_{1,i},y_{0,i},y_{1,i}:i\in {\Bbb
Z}_p)^{\langle\sigma_1\rangle}= K(\zeta)(x_i,y_i,X_i,Y_i:
i\in{\Bbb Z}_p)$ and the actions are given by
\begin{eqnarray*}
\sigma_2: &\hspace*{-4mm}& x_0\mapsto x_0,\ y_0\mapsto y_0,\ X_0\mapsto\zeta^2 X_0,\ Y_0\mapsto\zeta^{-2}Y_0, \\
          &\hspace*{-4mm}& x_i\mapsto\zeta x_i,\ y_i\mapsto\zeta^{-1}y_i,\ X_i\mapsto\zeta^{-1}X_i,\ Y_i\mapsto\zeta Y_i, \\
\sigma_3: &\hspace*{-4mm}& x_0\mapsto x_0,\ y_0\mapsto y_0,\ X_0\mapsto\zeta^2 X_0,\ Y_0\mapsto\zeta^{-2}Y_0, \\
          &\hspace*{-4mm}& x_i\mapsto\zeta x_i,\ y_i\mapsto\zeta^{-1}y_i,\ X_i\mapsto\zeta X_i,\ Y_i\mapsto\zeta^{-1}Y_i, \\
\sigma_4: &\hspace*{-4mm}& x_0\mapsto x_0x^p_1,\ y_0\mapsto
y_0y_1X_1,\
                           X_0\mapsto X_0X_1x_1,\ Y_0\mapsto Y_0Y_1y_1, \\
          &\hspace*{-4mm}& x_1\mapsto x_2\mapsto x_3\mapsto\cdots\mapsto x_{p-1}\mapsto(x_1x_2\cdots x_{p-1})^{-1}, \\
          &\hspace*{-4mm}& \mbox{Similarly for ~} y_1,y_2,\ldots,y_{p-1},X_1,\ldots,X_{p-1},Y_1,\ldots,Y_{p-1}. \\
\tau: &\hspace*{-4mm}& x_0\mapsto y^p_0x^{-1}_0,\ y_0\mapsto y_0,\ X_0\mapsto Y_0\mapsto X_0, \\
      &\hspace*{-4mm}& x_i\leftrightarrow y_i,\ X_i\leftrightarrow Y_i \mbox{ ~for~ } 1\le i\le p-1.
\end{eqnarray*}

\

Step 3. Define $\tilde{x}=x_0y^{-(p-1)/2}_0$,
$\tilde{y}=x^{-1}_0y^{(p+1)/2}_0$. Apply Theorem 3.1 again. It
suffices to consider $K(\zeta)(X_0,Y_0,x_i,y_i,X_i,Y_i: 1\le i\le
p-1)^{\langle\sigma_2,\sigma_3,\sigma_4,\tau\rangle}$.

We may apply Theorem 3.1 again to ``neglect" $X_0$ and $Y_0$. Thus
it suffices to consider $K(\zeta)(x_i,y_i,X_i,Y_i:1\le i\le
p-1)^{\langle\sigma_2,\sigma_3,\sigma_4,\tau\rangle}$.

\

Step 4. Define $u_0=x^p_1$, $v_0=x_1y_1$, $U_0=x_1X_1$,
$V_0=y_1Y_1$; and for $1\le i\le p-2$, define
$u_i=x^{-1}_ix_{i+1}$, $v_i=y^{-1}_i y_{i+1}$, $U_i=X^{-1}_i
X_{i+1}$, $V_i=Y^{-1}_i Y_{i+1}$. Then
$K(\zeta)(x_i,y_i,X_i,Y_i:1\le i\le p-1)^{\langle\sigma_2\rangle}=
K(\zeta)(u_i,v_i,U_i,V_i: 0\le i\le p-2)$. Note that the actions
of $\sigma_3$, $\sigma_4$, $\tau$ are given by
\begin{eqnarray*}
\sigma_3: &\hspace*{-4mm}& u_0\mapsto u_0,\ v_0\mapsto v_0,\ U_0\mapsto\zeta^2U_0,\ V_0\mapsto\zeta^{-2}U_0, \\
          &\hspace*{-4mm}& \mbox{All the other generators are fixed by ~} \sigma_3. \\
\sigma_4: &\hspace*{-4mm}& u_0\mapsto u_0u^p_1,\ v_0\mapsto
v_0v_1u_1,\
                           U_0\mapsto U_0U_1u_1,\ V_0\mapsto V_0V_1v_1, \\
          &\hspace*{-4mm}& u_1\mapsto u_2\mapsto\cdots\mapsto u_{p-2}\mapsto(u_0u^{p-1}_1u^{p-2}_2\cdots u^2_{p-2})^{-1}, \\
          &\hspace*{-4mm}& v_1\mapsto v_2\mapsto\cdots\mapsto v_{p-2}\mapsto u_0(v^p_0v^{p-1}_1\cdots v^2_{p-2})^{-1}, \\
          &\hspace*{-4mm}& U_1\mapsto U_2\mapsto\cdots\mapsto U_{p-2}\mapsto u_0(U^p_0U^{p-1}_1\cdots U^2_{p-2})^{-1}, \\
          &\hspace*{-4mm}& V_1\mapsto V_2\mapsto\cdots\mapsto V_{p-2}\mapsto u^{-1}_0v^p_0(V^p_0V^{p-1}_1\cdots V^2_{p-2})^{-1}, \\
\tau: &\hspace*{-4mm}& u_0\mapsto v^p_0u^{-1}_0,\ v_0\mapsto v_0,\ U_0\mapsto V_0\mapsto U_0, \\
      &\hspace*{-4mm}& u_i\leftrightarrow v_i,\ U_i\leftrightarrow V_i \mbox{ ~for ~} 1\le i\le p-2.
\end{eqnarray*}

\

Step 5. Define $R_0=U^p_0$, $S_0=U_0V_0$; and for $1\le i\le p-2$,
define $R_i=U_i$, $S_i=V_i$. Then $K(\zeta)(u_i,v_i,U_i,V_i: 0\le
i \le p-2)^{\langle \sigma_3\rangle}= K(\zeta)(u_i,v_i,R_i,S_i:
0\le i\le p-2)$. We will write the actions of $\sigma_4$ and
$\tau$ on $R_i$, $S_i$ as follows.
\begin{eqnarray*}
\sigma_4: &\hspace*{-4mm}& R_0\mapsto R_0R^p_1u^p_1,\ S_0\mapsto S_0S_1R_1u_1v_1, \\
          &\hspace*{-4mm}& R_1\mapsto R_2\mapsto\cdots\mapsto R_{p-2}\mapsto u_0(R_0R^{p-1}_1R^{p-2}_2\cdots R^2_{p-2})^{-1}, \\
          &\hspace*{-4mm}& S_1\mapsto S_2\mapsto\cdots\mapsto S_{p-2}\mapsto u^{-1}_0v^p_0R_0(S^p_0S^{p-1}_1\cdots S^2_{p-2})^{-1}, \\
\tau: &\hspace*{-4mm}& R_0\mapsto S^p_0R^{-1}_0,\ S_0\mapsto S_0,\
R_i\leftrightarrow S_i\mbox{ ~for~ } 1\le i\le p-2.
\end{eqnarray*}

\

Step 6. Imitate the change of variables in Step 6 of Case 1. We
define $u_{p-1}=(u_0 u^{p-1}_1 u^{p-2}_2 \cdots u^2_{p-2})^{-1}$,
$R_{p-1}=u_0(R_0 R^{p-1}_1 R^{p-2}_2 \cdots R^2_{p-2})^{-1}$; and
for $1\le i\le p-2$, define $w_i=v_0 v_1 \cdots v_i u_1 u_2 \cdots
u_i$, $T_i=S_0 S_1 \cdots S_i R_1 R_2 \cdots R_i$; define
$w_{p-1}=(v^{p-1}_0$ $v^{p-2}_1 \cdots v_{p-2} u^{p-2}_1 u^{p-3}_2
\cdots u_{p-2})^{-1}$, $T_{p-1}=u_1v_1v^p_0(S^{p-1}_0 S^{p-2}_1
\cdots S_{p-2} R^{p-2}_1 R^{p-3}_2 \cdots R_{p-2})^{-1}$.

We find that $K(u_i,v_i,R_i,S_i: 0\le i\le
p-2)=K(u_i,w_i,R_i,T_i:1\le i\le p-1)$ and
\begin{eqnarray*}
\sigma_4: &\hspace*{-4mm}& u_1\mapsto u_2\mapsto\cdots\mapsto u_{p-1}\mapsto(u_1u_2\cdots u_{p-1})^{-1}, \\
          &\hspace*{-4mm}& w_1\mapsto w_2\mapsto\cdots\mapsto w_{p-1}\mapsto(w_1w_2\cdots w_{p-1})^{-1}, \\
          &\hspace*{-4mm}& R_1\mapsto R_2\mapsto\cdots\mapsto R_{p-1}\mapsto(R_1R_2\cdots R_{p-1})^{-1}, \\
          &\hspace*{-4mm}& T_1\mapsto T_2\mapsto\cdots\mapsto T_{p-1}\mapsto(T_1T_2\cdots T_{p-1})^{-1}, \\
\tau: &\hspace*{-4mm}& w_i\mapsto w_i,\ T_i\mapsto T_i,\
                       u_i\mapsto w_i(u_iw_{i-1})^{-1},\ R_i\mapsto T_i(R_i T_{i-1})^{-1}
\end{eqnarray*}
where we write $w_0=v_0$, $T_0=S_0$ for convenience.

\

Step 7. The multiplicative action in Step 6 can be formulated as
follows.

Let $\pi=\langle \tau,\sigma_4\rangle$ and define a $\pi$-lattice
as in Step 7 of Case 1 (but $\sigma_3$ should be replaced by
$\sigma_4$ in the present situation). Then
$K(\zeta)(u_i,w_i,R_i,T_i:1\le i\le p-1)^\pi =K(\zeta)(M\oplus
M)^\pi$ where $M$ is the same lattice in Step 7 of Case 1.

The structure of $M$ has been determined in Step 7 of Case 1. Thus
$M\oplus M\simeq \Lambda \oplus\Lambda$ where $\Lambda\simeq {\Bbb
Z}[\pi_1] \otimes_{\Bbb Z} {\Bbb Z}[\pi_2]/\Phi_p(\sigma_4)$ where
$\pi_1=\langle\tau\rangle$ and $\pi_2=\langle\sigma_4\rangle$. It
follows that we can find elements $Y_1,\ldots,Y_{p-1}$,
$Z_1,\ldots,Z_{p-1}$, $W_1,\ldots,W_{p-1}$, $Q_1,\ldots,Q_{p-1}$
in the field $K(\zeta)(u_i,w_i,R_i,T_i:1\le i\le p-1)$ so that
$K(\zeta)(u_i,w_i,R_i,T_i:1\le i\le
p-1)=K(\zeta)(Y_1,\ldots,Y_{p-1},Z_1,\ldots,Z_{p-1},W_1,\ldots,W_{p-1}$,
$Q_1,\ldots,Q_{p-1})$ and the actions of $\sigma_4$ and $\tau$ are
given by
\begin{eqnarray*}
\sigma_4: &\hspace*{-4mm}& Y_1\mapsto \cdots\mapsto Y_{p-1}\mapsto(Y_1Y_2\cdots Y_{p-1})^{-1}, \\
          &\hspace*{-4mm}& \mbox{Similarly for ~} Z_1,\ldots,Z_{p-1},\ W_1,\ldots,W_{p-1},\ Q_1,\ldots,Q_{p-1}; \\
\tau: &\hspace*{-4mm}& Y_i\leftrightarrow Z_i,\ W_i\leftrightarrow
Q_i.
\end{eqnarray*}

\

Step 8. We can linearize the actions of $\sigma_4$ and $\tau$ by
the same method in Step 9 of Case 1. Apply Theorem 3.1. We may
``neglect" the roles of $W_i$, $Q_i$. Thus the present situation
is the same as in Step 8 of Case 1.\hfill$\square$

\newpage
\renewcommand{\refname}{\centering{References}}

\end{document}